\documentclass[a4paper]{amsart}
\usepackage[margin=3.5cm]{geometry}
\usepackage{amsmath,amssymb,amsthm,amscd,amsfonts,comment,mathtools}

\usepackage{tikz}
\usetikzlibrary{cd}
\usetikzlibrary{fit, patterns}
\usetikzlibrary{arrows,decorations.markings,plotmarks}

\newtheorem{theorem}{Theorem}[section]

\theoremstyle{definition}
\newtheorem{definition}[theorem]{Definition}

\theoremstyle{remark}
\newtheorem{remark}[theorem]{Remark}

\newtheorem{example}[theorem]{Example}

\DeclareMathOperator{\Res}{Res} 
\DeclareMathOperator{\RE}{R}   
\DeclareMathOperator{\Q}{\mathbf Q}
\DeclareMathOperator{\Z}{\mathbf Z}
\DeclareMathOperator{\E}{\mathcal E}
\DeclareMathOperator{\C}{\mathbf C}
\DeclareMathOperator{\R}{\mathbf R}
\DeclareMathOperator{\GL}{GL}

\DeclareMathOperator{\codim}{codim}
\DeclareMathOperator{\Tp}{Tp}
\DeclareMathOperator{\KTp}{KTp}
\DeclareMathOperator{\KTs}{KTs}
\DeclareMathOperator{\Ts}{Ts}

\DeclareMathOperator{\Ell}{Ell}
\DeclareMathOperator{\csm}{csm}
\DeclareMathOperator{\ssm}{ssm}
\DeclareMathOperator{\ssmTp}{ssm-\!Tp}
\DeclareMathOperator{\mSTp}{ms-\!Tp}
\DeclareMathOperator{\h}{\hbar}
\DeclareMathOperator{\MC}{mc}
\DeclareMathOperator{\MS}{ms}
\DeclareMathOperator{\Var}{Var}

\def\pt{pt}

\title{Thom polynomials. A primer}
\author{R. Rim\'anyi}

\begin{document}

\begin{abstract}
The Thom polynomial of a singularity $\eta$ expresses the cohomology class of the $\eta$-singularity locus of a map in terms of the map's simple invariants. In this informal survey---based on two lectures given at the Isaac Newton Institute in 2024---we explore various Thom polynomial concepts with examples.
\end{abstract}

\maketitle

\section{Introduction}\label{sec:intro}
\begin{figure}
\[
\begin{tikzpicture}[scale=.5]
   \draw (0,0) ellipse (5 and 3);
   \draw (-3,.9) to[out=-50,in=180] (0,-.7) to[out=0,in=-130] (3,.9);
   \draw (-2.1,0) to[out=33,in=180] (0,.7) to[out=0,in=147] (2.1,0);
   \draw[dashed] (-2.1,0) to[out=213,in=50] (-3,-0.9);
   \draw[dashed] (-3,-0.9) to[out=51,in=-90] (-2.6,0) to[out=90,in=-51] (-3,0.9);
   \draw[dashed] (2.1,0) to[out=-33,in=130] (3,-0.9);
   \draw[dashed] (3,-0.9) to[out=129,in=-90] (2.6,0) to[out=90,in=241] (3,0.9);
   \draw[blue] (8,2) node {$A_0$-point};
   \draw[blue,->] (6.5,2) to[out=160,in=20] (1,1.5);
   \draw[blue] (8,0) node {$A_1$-point};
   \draw[blue,->] (6.5,0) to[out=160,in=20] (5.1,0);
   \draw[blue] (8,-2) node {$A_2$-point};
   \draw[blue,->] (6.5,-2) to[out=200,in=-40] (3,-1.1);
\end{tikzpicture}
\]
\caption{A map of the torus to the plane. The map stratifies the source (the torus) to $A_0$-, $A_1$-, and $A_2$-points. The $A_0$-points are dense open, the $A_1$-points form a circle and four open intervals, there are four $A_2$-points.}
\label{fig:torus}
\end{figure}
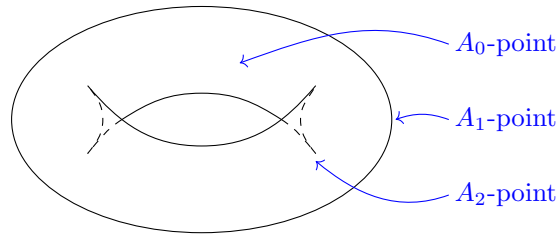

Let us consider the usual picture of the torus $T=S^1\times S^1$ in Figure \ref{fig:torus}. Our goal with this illustration is not only to represent the torus $ T $ but also to show a  (smooth) map $ f: T \to \R^2 $, where $\R^2$ represents the flat piece of paper where the torus image appears. For clarity, we have sketched the torus as transparent.

As a consequence of the Inverse Function Theorem, around most points of $T$ this map can be locally described by $(x,y)\mapsto (x,y)$. That is, local coordinate systems can be chosen around such a point on the torus, and around its image in $\R^2$, so that the map $f$ in these coordinate systems is the identity map. We call such points $A_0$-points of the map $f$.

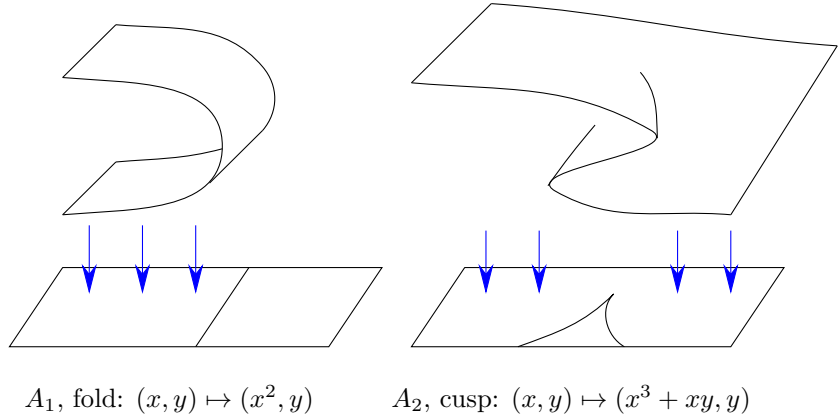
\begin{figure}
\[
\begin{tikzpicture}[scale=0.7]
    \draw (-1,-1) to[out=5,in=-90] (2,.3) to[out=90,in=-5] (-1,1.6);
    \draw (0,0) to[out=5,in=195] (2,0.25);
    \draw (2.76,.6) to[out=50,in=-50] (2.76,1.8) to[out=130,in=-5] (0,2.6);
    \draw (-1,-1) to (0,0); 
    \draw (-1,1.6) to (0,2.6);
    \draw (1.76,-.4) to (2.76,.6);
    \draw (-1,-2) to (5,-2) to (4,-3.5) to (-2,-3.5) to (-1,-2);
    \draw (2.5,-2) to (1.5,-3.5);
\draw[-{Stealth[length=4mm, width=2mm]},blue] (-.5,-1.2) to (-.5,-2.5);
\draw[-{Stealth[length=4mm, width=2mm]},blue] (.5,-1.2) to (.5,-2.5);
\draw[-{Stealth[length=4mm, width=2mm]},blue] (1.5,-1.2) to (1.5,-2.5);
\draw (1,-4.5) node {$A_1$, fold: $(x,y)\mapsto (x^2,y)$};
\end{tikzpicture}
\begin{tikzpicture}[scale=.7]
\draw [black](-2,1.5) to [out=-5,in=150] (2.5,.6) to [out=-30, in=150] (.7,-.55) to [out=-30,in=175] (4,-1);
\draw [black](-2,1.5) to (-.5,3);
\draw [black](-.5,3) to[out=-4,in=176] (6,2.2);
\draw [black] (6,2.2) to (4,-1);
\draw [black] (2.62,0.45) to[out=92,in=-50] (2.3,1.7); 
\draw [black] (.565,-0.45) to[out=55,in=-130] (1.45,.7);
\draw [black] (-1,-2) to (5,-2) to (4,-3.5) to (-2,-3.5) to (-1,-2);
\draw (0,-3.5) to[out=20,in=-135] (1.8,-2.5) to[out=-135,in=150] (2,-3.5);
\draw[-{Stealth[length=4mm, width=2mm]},blue] (-.6,-1.3) to (-.6,-2.5);
\draw[-{Stealth[length=4mm, width=2mm]},blue] (0.4,-1.3) to (0.4,-2.5);
\draw[-{Stealth[length=4mm, width=2mm]},blue] (3,-1.3) to (3,-2.5);
\draw[-{Stealth[length=4mm, width=2mm]},blue] (4,-1.3) to (4,-2.5);
\draw (1,-4.5) node {$A_2$, cusp: $(x,y)\mapsto (x^3+xy,y)$};
\end{tikzpicture}
\]
\caption{The of $A_1$ and $A_2$ singularities $(\R^2,0)\to (\R^2,0)$.}  
\label{fig:A1A2}
\end{figure}

Around certain other points of the torus the map can be locally written as $(x,y)\mapsto (x^2,y)$---see Figure \ref{fig:A1A2} (left). We call these points the $A_1$-points, or {\em fold} points of $f$. This set of points consists of a circle and four open intervals. Lastly, at four special points on the torus the map can be written locally as $(x,y)\mapsto (x^3+xy,y)$, see Figure \ref{fig:A1A2} (right).  These points we call the $A_2$-points, or {\em cusp} points of the map $f$.

What we described is an example of the stratification of the source manifold of a map $f:M^m\to N^n$ to strata called singularity loci. 

\medskip

Thom polynomials express the cohomology class represented by the closure of such singularity loci in terms of simple invariants of the spaces and the map. More precisely, for each singularity $\eta,$ there is an associated multivariate polynomial $\Tp(\eta),$ known as the Thom polynomial. By substituting the characteristic classes of the source manifold and the target manifold (pulled back to the source by the map) into this polynomial, we obtain the cohomology class represented by the closure of the $\eta$-singularity locus.

An important point is that $\Tp(\eta)$ does not depend on the specific situation; it depends only on the singularity. This property is sometimes emphasized by calling $\Tp(\eta)$ a ``universal polynomial." For example, the Thom polynomial of the cusp singularity discussed above is $ a_2 + a_1 b_1 + b_1^2 + b_2, $ where $ a_i $ are the Stiefel-Whitney classes of the domain and $ b_i $ are the Stiefel-Whitney classes of the target surface (pulled back by $ f^* $). If the target surface is the plane or the sphere (where $ b_i = 0 $), we find that the parity of the cusp points of $ f: M^2 \to \R^2 $ is the same as the parity of the Euler characteristic of $ M^2 $ (that is, $ a_2 $). In particular, a generic map of $T^2$ to $\R^2$ has an even number of cusps, while a generic map of the  real projective plane to $\R^2$ has an odd number of cusps.

While the previous example was over the reals, most of our discussions will be over the complex numbers, and in higher dimensions. Consequently, Chern classes, rather than Stiefel-Whitney classes, will be our characteristic classes.

\medskip

\noindent {\bf Plan of the paper.} We will first discuss the relevant concepts from singularity theory. Then, we will review various versions of Thom polynomials. Throughout this survey, we will emphasize examples over precisely phrased theorems. Our treatment of the subject is both informal and incomplete (see Section \ref{sec:whatelse}) and is shaped by the author's knowledge and mathematical preferences. For a more formal and comprehensive survey, readers are advised to consult \cite{OhmotoSurvey}.

\medskip

\noindent{\bf Acknowledgments.}
This article is the written version of the lectures given at the {\em Singularity theory and hyperbolicity} workshop at the Isaac Newton Institute for Mathematical Sciences in the Spring of 2024. The author is grateful to the organizers, as well as the Newton Institute for their hospitality. 
This work was partially supported by NSF grants  2152309, 2200867 and a grant from the Heilbronn Institute. Special thanks to G.~B\'erczi, L.~Feh\'er, and T.~Ohmoto for helpful discussions during the preparation of this paper.

\section{Singularities of maps}

A general reference of singularity theory is \cite{AVGL}.

\subsection{Singularities---first attempt} \label{sec:Asing}
Let us discuss what we should mean by a ``singularity'' $\eta$ so that the construction of {\em $\eta$-points of a map} in the last section makes sense.

Let $\E(m,n)$ denote the vector space of holomorphich map {\em germs} from $(\C^m,0)$ to $(\C^n,0)$. We will assume that $m\leq n$ and set $\ell=n-m\geq 0$. Germs are equivalence classes of maps $U\to \C^n$ defined in a neighborhood $U$ of $0\in\C^m$, satisfying $0\mapsto 0$. Two such maps represent the same germ if they agree in a neighborhood of $0\in \C^m$. 
A global map $M^m\to N^n$, at point $p\in M$, defines an element of $\E(m,n)$ only if coordinate charts are fixed around $p$ and $f(p)$. Since there are no a priori choices of charts on a manifold, we need to build this ambiguity in our definition of ``singularity''.

Let $\E^o(m)\subset \E(m,m)$ denote the subset of invertible germs, and observe that it is a group by composition. Moreover, we have the action of $\E^o(m)\times \E^o(n)$ on $\E(m,n)$ by
\begin{equation}\label{def:RLaction}
    (\alpha,\beta)\cdot f=\beta \circ f \circ \alpha^{-1},
\end{equation}
and it encodes exactly the ambiguity of the choice of charts in the source and in the target. 

\begin{definition}\label{def:Asing}
 An $\mathcal A$-singularity (or ``right-left singularity'') $\eta$ is an orbit of the action \eqref{def:RLaction}.
\end{definition}

With this definition the construction of Section \ref{sec:intro} makes sense: let $\eta$ be an $\mathcal A$-singularity from $m$ to $n$ dimensions, and let $f:M^m\to N^n$ be a map between complex manifolds. Define the {\em singularity locus}
\[
\eta(f)=\{
p\in M : \text{the germ of $f$ at $p$ belongs to } \eta\}. 
\]

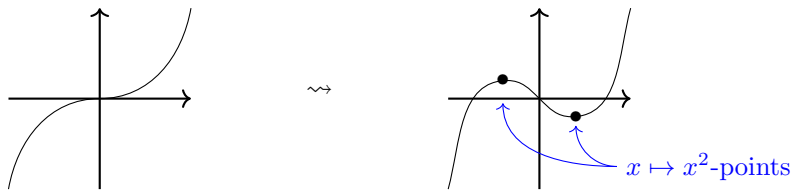
\begin{figure}
\[
    \begin{tikzpicture}[scale=.6,baseline=1.1]
        \draw[thick,->] (-2,0) -- (2,0);
        \draw[thick,->] (0,-2) -- (0,2);
        \draw (-2,-2) to[out=80,in=180] (0,0) to[out=0,in=260] (2,2);
    \end{tikzpicture}
    \qquad\qquad\leadsto\qquad\qquad
    \begin{tikzpicture}[scale=.6,baseline=1.1]
        \draw[thick,->] (-2,0) -- (2,0);
        \draw[thick,->] (0,-2) -- (0,2);
        \draw (-2,-2) to[out=75,in=180] (-.7,.4) to [out=0,in=135] (0,0) to[out=-45,in=180] (0.7,-0.4) to
        [out=0,in=255] (2,2);
        \draw (-.8,.4) node {$\bullet$};
        \draw (.8,-.4) node {$\bullet$};
        \draw[blue] (3.7,-1.5) node {$x\mapsto x^2$-points};
        \draw[blue,->] (1.7,-1.5) to[out=180,in=-90] (-0.8,-0.1);
        \draw[blue,->] (1.7,-1.5) to[out=180,in=-90] (0.8,-0.6);
    \end{tikzpicture}
    \]
    \caption{The singularity $x\mapsto x^3$ is not stable, it disappears at the perturbation $x\mapsto x^3-\varepsilon x$. (For $\varepsilon<0$ the appearing $x\mapsto x^2$-points are complex.)}
    \label{fig:x3}
\end{figure}

While this definition is mathematically sound, it presents a few practical issues. Firstly, there are an overwhelming number of $\mathcal{A}$---singularities--—trust me on this.
Another issue is that some $\mathcal{A}$-singularities are irrelevant for us. For example, when $m=n=1$, the germ of $x\mapsto x^3$ does not appear in a generic map. This is because if we perturb $x\mapsto x^3$, it splits into two $x\mapsto x^2$ singularities, as illustrated in Figure \ref{fig:x3}. Our term for this phenomenon is that the $\mathcal A$-singularity $x\mapsto x^3$ is not {\em stable}. For our purposes non-stable singularities will be irrelevant. 

The $\mathcal A$-singularities $x\mapsto x$, $x\mapsto x^2$, as well as the three singularities (from two to two dimensions) appearing in Section \ref{sec:intro} are stable.

\begin{remark}
\label{rem:jets}
    According to our definition, $\mathcal{A}$-singularities are orbits of an infinite-dimensional group acting on an infinite-dimensional vector space, However, for practical purposes, we can truncate the elements of the group and the vector space at a given order $N$. The technical term for this truncation is the use of {\em $N$-jets} instead of germs. By doing so, the group and the vector space become finite-dimensional. Although not every $\mathcal{A}$-singularity can be detected in a jet space for large $N$, the ones that cannot be detected are irrelevant for the theory of Thom polynomials.
\end{remark}

\subsection{Singularities parametrized by algebras}
Experience shows that the right notion of singularities is achieved if we ``glue together'' the $\mathcal A$-singularities with the same local algebra.

\begin{definition}
The local algebra $Q_f$ of a germ 
\[
f:(x_1,\ldots,x_m)\mapsto (f_1,\ldots,f_n)
\]
in $\E(m,n)$ is defined to be 
\[
\C[\![ x_1,\ldots,x_m ]\!]/(f_1,\ldots,f_n).
\]
\end{definition}

We will assume in the whole paper that $Q_f$ is finite dimensional. Again, just like the germs that are not defined by polynomials, the ones for which $Q_f$ is infinite dimensional are not relevant in Thom polynomial theory.

\begin{example}
 For $f_1: x \mapsto x^2$, $f_2:x\mapsto x^3$, $f_3:x\mapsto x^2+x^3$ we have
 \[Q_{f_1}\cong Q_{f_3}\cong\C[x]/(x^2), \qquad
 Q_{f_2}\cong\C[x]/(x^3).\]
 For $f_4:(x,y)\mapsto (x,y)$, $f_5:(x,y)\mapsto (x^2,y)$, $f_6:(x,y)\mapsto (x^3+xy,y)$, ${f_7}:(x,y)\mapsto (x^3,y)$ we have 
 \[
 Q_{f_4}\cong\C, \qquad 
 Q_{f_5}\cong\C[x]/(x^2), \qquad
 Q_{f_6}\cong Q_{f_7} \cong \C[x]/(x^3).
 \]
 Observe that $f_6$ and $f_7$ do not belong to the same $\mathcal A$-singularity, but their local algebras are isomorphic. In fact, $f_6$ is stable while $f_7$ is not.
 \end{example}

 \begin{example} \label{ex:I22}
 The local algebra of any of the germs $f_8:(x,y)\mapsto (xy,x^2+y^2)$, $f_9:(x,y)\mapsto (x^2, y^2)$, $f_{10}:(x,y,u,v)\mapsto (xy,x^2+y^2+ux+vy,u,v)$ is isomorphic to 
$\C[x,y]/(xy,x^2+y^2)$. We denote this algebra by $I_{22}$. The germs $f_8$ and $f_9$ are not stable, while $f_{10}$ is stable. In fact there are no stable $\mathcal A$-singularities $(\C^2,0)\to (\C^2,0)$ with local algebra $I_{22}$. There are stable singularities $(\C^m,0)\to (\C^n,0)$ with local algebra $I_{22}$ if and only if $\ell\geq0$ and $m\geq 3\ell+4$ (recall $\ell=n-m$).
\end{example}

\begin{definition}
\label{def:contact}
    Let $Q$ be a finite dimensional local algebra, and $m\leq n$ non-negative integers. Define the ``contact singularity'' $\eta_{Q}=\eta_Q(m,n)\subset \E(m,n)$ to be the set of germs whose local algebra is isomorphic to $Q$.
\end{definition}

This definition may sound arbitrary, but in fact there is a geometric story behind it. Namely, one can define a group $\mathcal K=\mathcal K(m,n)$, the ``contact group'' acting on $\E(m,n)$ that contains $\E^o(m)\times \E^o(n)$. Graphs of $\mathcal K$-equivalent germs have the same {\em contact} with $\C^m \times \{0\}\subset \C^m\times \C^n$---hence the name. It is a theorem of Mather that belonging to the same $\mathcal K$-orbit is equivalent to having isomorphic local algebras---hence Definition \ref{def:contact} is in fact geometric!

A second compelling case for Definition \ref{def:contact} is another classical theorem: given a $Q$, for large enough $\ell$ and for large enough $m$ compared to $\ell$, the set $\Sigma_Q$ contains a dense-open {\em stable} $\mathcal A$-singularity. That is, closures of $\eta_Q$ and closures of stable $\mathcal A$-singularities are the same for large $\ell$ and $m$. 

\begin{remark}
    How large do $\ell$ and $m$ need to be for the mentioned theorem? The thresholds are calculable invariants of the local algebra $Q$---the interested reader will find the relevant theorems and algorithms when studying ``miniversal unfoldings'', as well as ``genotype'' and ``prototype'' germs for the algebra $Q$. One example is given in Example \ref{ex:I22}, here is another one: for $Q=\C[x,y]/(x^2,xy,y^2)$ we need $\ell \geq 1$ and $m\geq 2\ell+4$.
\end{remark}

\subsection{The zoo of contact singularities for a given $\ell$}
\begin{figure}
\[
\tikz[overlay]{
\filldraw[fill=yellow!30,draw=black!100]  (1.3,-6.3) [rounded corners]to
  (1.6,4.3)  to
  (3.2,4.3)  [sharp corners] to 
  (3.5,-6.3);
\filldraw[fill=yellow!35,draw=black!100]  (1.4,-6.3) [rounded corners]to
  (1.6,3.1)  to
  (3.2,3.1)  [sharp corners] to 
  (3.4,-6.3);
\filldraw[fill=yellow!40,draw=black!100]  (1.4,-6.3) [rounded corners]to
  (1.6,1.9)  to
  (3.2,1.9)  [sharp corners] to 
  (3.4,-6.3);
\filldraw[fill=yellow!45,draw=black!100]  (1.4,-6.3) [rounded corners]to
  (1.6,.7)  to
  (3.2,.7)  [sharp corners] to 
  (3.4,-6.3);
\filldraw[fill=yellow!50,draw=black!100]  (1.4,-6.3) [rounded corners]to
  (1.6,-.5)  to
  (3.2,-.5)  [sharp corners] to 
  (3.4,-6.3);
\filldraw[fill=yellow!55,draw=black!100]  (1.4,-6.3) [rounded corners]to
  (1.6,-1.7)  to
  (3.2,-1.7)  [sharp corners] to 
  (3.4,-6.3);
\filldraw[fill=yellow!60,draw=black!100]  (1.4,-6.3) [rounded corners]to
  (1.6,-2.9)  to
  (3.2,-2.9)  [sharp corners] to 
  (3.4,-6.3);
\filldraw[fill=yellow!65,draw=black!100]  (1.4,-6.3) [rounded corners]to
  (1.6,-4.1)  to
  (3.2,-4.1)  [sharp corners] to 
  (3.4,-6.3);  
\filldraw[fill=yellow!70,draw=black!100]  (1.4,-6.3) [rounded corners]to
  (1.6,-5.3)  to
  (3.2,-5.3)  [sharp corners] to 
  (3.4,-6.3);
\filldraw[fill=blue!20,draw=black!60]     (3.5,-6.3) [rounded corners]to
  (3.6,.7)  to
  (4.5,.7) to
  (8,-2)  to 
  (11,-3) [sharp corners] to 
  (11.5,-6.3);
\filldraw[fill=blue!30,draw=black!60]     (8.4,-6.3) [rounded corners]to
  (8.4,-2.5)  to
  (10.7,-3.1) [sharp corners] to
  (11.3,-6.3);
\filldraw[fill=red!20,draw=black!60]     (11.6,-6.3) [rounded corners]to
  (11.6,-5.5)  to
  (12.5,-5.5) [sharp corners] to
  (13,-6.3);
\draw[blue] (4.5,4.5) node {$\Sigma^1$};
\draw[->, blue] (4.2,4.5) to[out=170,in=40] (3.2,4.2);

\draw[blue] (4.8,3.5) node {$\Sigma^{11}$};
\draw[->, blue] (4.4,3.5) to[out=170,in=80] (2.9,3.1);

\draw[blue] (5.1,2.5) node {$\Sigma^{111}$};
\draw[->, blue] (4.7,2.5) to[out=170,in=80] (2.9,1.9);

\draw[blue] (5.1,2) node {$\cdots$};

\draw[blue] (6.5, 1) node {$\Sigma^2$};
\draw[->, blue] (6.2,1) to[out=190,in=60] (5.05,.39);

\draw[blue] (9.4, -1.1) node {$\Sigma^{21}$};
\draw[->, blue] (9.1,-1.3) to[out=-130,in=80] (8.7,-2.55);

\draw[blue] (12, -2.5) node {$\Sigma^{3}$};
\draw[->, blue] (12,-2.8) to[out=-100,in=90] (11.9,-5.4);

\draw[red] (12,0) node {$\lambda$: moduli!};
\draw[->, red] (12,-.3) to[out=-60,in=70] (12.4,-6.1);
}
\begin{tikzcd}
\codim \\
0 & A_0 \\
1 & A_1 \\ 
2 & A_2 \\
3 & A_3 \\
4 & A_4 & I_{22} \\
5 & A_5 & I_{23} \\
6 & A_6 & I_{24} & I_{33} \\
7 & A_7 & I_{25} & I_{34} & & (x^2,y^3) \\
8 & A_8 & I_{26} & I_{35} & I_{44} & (x^2+y^3,xy^2) \\
9 & A_9 & & \cdots & & \cdots  & C_\lambda
\end{tikzcd}
\]
\caption{The zoo of small codimension contact singularities for $\ell=0$. Notation: $A_m=\C[x]/(x^{m+1})$, $I_{ab}=\C[x,y]/(xy,x^a+y^b)$, 
$(x^2,y^3)$ stands for $\C[x,y]/(x^2,y^3)$, $(x^2+y^3,xy^2)$ stands for $\C[x,y]/(x^2+y^3,xy^2)$.
$C_\lambda=\C[x,y,z]/(x^2-\lambda yz, y^2-\lambda xz, z^2-\lambda xy)$, that is a codimension 10 singularity for each concrete $\lambda$ with $\lambda(\lambda^3-1)(8\lambda^3+1)\not=0$---their union is 9 codimensional. The $\Sigma^I$-decorations of the figure are discussed in Section~\ref{sec:ThomBoardman}.
}
\label{Fig:l0}
\end{figure}
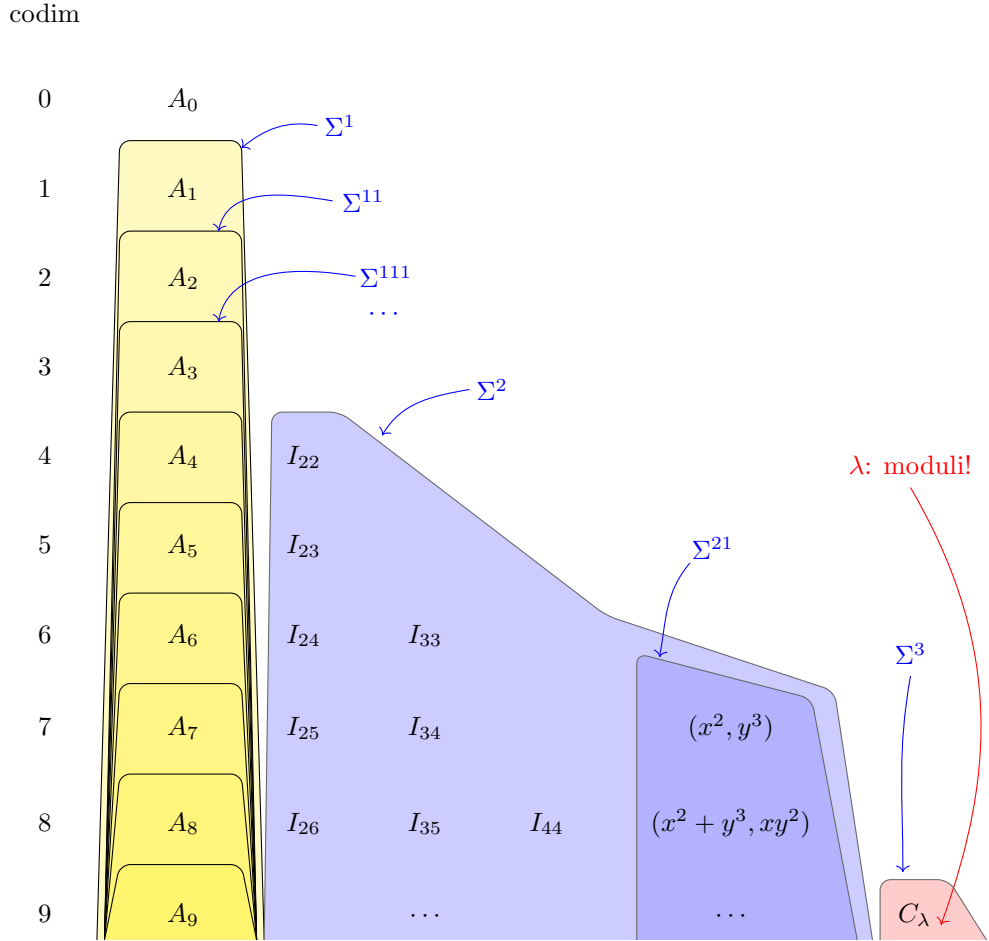

Luckily the hierarchy of contact singularities in $\E(m,n)$ are not independent for any $m$ and $n$, but in fact essentially depend on $\ell=n-m$ only. For example, for $\ell=0$ the contact singularities of codimension at most 8 are in Figure \ref{Fig:l0}. For a concrete $m$ only the codimension $\leq m$ ones appear for a stable map. For example, for maps $M^2\to N^2$ we see the three possible singularities, $A_0$, $A_1$, $A_2$ of Section \ref{sec:intro}.

It is important to point out a new phenomenon: in codimension 9 the hierarchy starts to contain moduli. The algebras $\C[x,y,z]/(x^2-\lambda yz, y^2-\lambda xz, z^2-\lambda xy)$ are typically not isomorphic for different $\lambda$'s.

The task of classifying all algebras, even the ones that occur as local algebras of singularities for a given $\ell$, is of course hopeless. Yet, in various special cases the classification is known, see eg. \cite{dPW}

It is natural to study the relation of such hierarchies for different $\ell$'s. The first fact to keep in mind is that as we increase $\ell$ new algebras appear. For example the algebra $III_{22}=\C[x,y]/(x^2,xy,y^2)$ only appears for $\ell\geq 1$. (The reason is that its minimal presentation has {\em one} more relations than generators.) 
Another fact is that the relative position of algebras may change considerably as we increase $\ell$. Let us illustrate this with pointing out that 
\[
\codim( A_6 ) =6\ell + 6,
\qquad
\codim(\C[x,y]/(x^2,y^3))=5\ell +7,
\]
hence, for $\ell=0$ the singularity $A_6$ is ``above'' the other one, but for large $\ell$ it is way ``below'' it.

When $\ell$ is large enough so that $\eta_Q$ exists in $\E(m,m+\ell)$ then the formula for its codimension is 
\[
\codim(\eta_Q)=(\dim_{\C}(Q)-1)\ell+b(Q)
\]
where $b(Q)$ is also a computable algebraic invariant of $Q$. 
Hence as $\ell$ is large, the top of the hierarchy of contact singularities starts with the small dimensional algebras, see Figure~\ref{fig:largel}.

\begin{figure} 
\begin{center}
\begin{tabular}{l|clclll|c}
 & &   & & & &  & appears for: \\
$\dim(Q)=1$ & &   $\codim(A_0)$ & $=$ & $0\ell +0$ &&& $\ell\geq 0$
    \\
\hline
$\dim(Q)=2$ & &    $\codim(A_1)$ & $=$ & $1\ell +1$ &&& $\ell\geq 0$\\
\hline
$\dim(Q)=3$ & &    $\codim(A_2)$ & $=$ & $2\ell+2$ &&& $\ell\geq 0$\\
& &    $\codim(III_{22})$ & $=$ & $2\ell +4$ &&& $\ell\geq 1$ \\
\hline
$\dim(Q)=4$ & &    $\codim(A_3)$ & $=$ & $3\ell +3$ &&& $\ell\geq 0$\\
& &    $\codim(I_{22})$ & $=$ & $3\ell +4$ &&& $\ell\geq 0$\\
& &   $\codim(III_{23})$ & $=$ & $3\ell+5$ &&& $\ell\geq 1$\\
& &    $\codim(\C[x,y,z]/(x,y,z)^2)$ & $=$ & $3\ell+9$ &&& $\ell\geq 6$\\
\hline
$\cdots$ & & $\cdots$ & & & & &
\end{tabular}
\end{center}
\caption{The top of the hierarchy of contact singularities for large $\ell$.}\label{fig:largel}
\end{figure}

\subsection{Thom-Boardman stratification}
\label{sec:ThomBoardman}
We started our journey to singularities with $\mathcal A$-singularities, and mentioned that there are too many of them. Hence, we glued some together and obtained contact singularities, and found that there are essentially as many of them as finite dimensional local algebras. One may argue that they are still too many. Indeed, there is a remarkable even rougher stratification: the Thom-Boardman singularities.

Thom-Boardman singularities of order 1 are 
\[
\Sigma^i=\{f \in \E(m,n) : \dim\ker f'(0) = i \}.
\]
Thom-Boardman singularities of order 2 are {intuitively} 
defined by 
\[
\Sigma^{ij}=\{f \in \E(m,n) : 
f\in \Sigma^i, f|_{\Sigma^i(f)} \in \Sigma^j \}.
\]
We advice the reader to verify that both the fold
$(x,y)\mapsto (x^2,y)$ and the 
cusp $(x,y)\mapsto (x^3+xy,y)$ belong to $\Sigma^1$, but the first one is $\Sigma^{10}$ the latter one is $\Sigma^{11}$.

The intuitive definition of higher order Thom-Boardman singularities is analogous. In Figure~\ref{Fig:l0} we indicated the Thom-Boardman types. 

There is a simple formula for the  codimension of Thom-Boardman singularities. For order~1 it is: $\codim(\Sigma^i)=i(i+\ell)$. The story of what contact singularities occupy ``most'' of a Thom-Boardman singularity has surprises. In the $\Sigma^2$ Thom-Boardman type $I_{22}$ is dense-open for $\ell=0$ but $III_{22}$ is dense-open for $\ell\geq 1$ (see the table above). A dense-open part of the $\Sigma^3$ Thom-Boardman type is a moduli of contact singularities for $\ell=0$ (see Figure \ref{Fig:l0}) but has a dense open contact singularity for $\ell\geq 6$ (the last line in Figure \ref{fig:largel}). For any multiindex $I$, for large enough $\ell$ the Thom-Boardman type $\Sigma^I$ contains a dense-open contact singularity.

\section{Thom polynomials}

In the 1950s, René Thom made significant contributions to both the local theory of singularities and the global theory of cobordisms and characteristic classes. He also bridged these seemingly distant areas, developing what is now known as the theory of Thom polynomials~ \cite{thomoriginal}.

\subsection{Definitions}

Let $\eta\subset \E(m,n)$ be a singularity. For this now we mean a (locally closed) subvariety invariant under at least the reparametrization group $\E^o(m) \times \E^o(n)$---for example an $\mathcal A$-singularity, a contact singularity, or a Thom-Boardman singularity. Let us assume it has pure (complex) codimension $d$.

The Thom polynomial of $\eta$
\[
\Tp(\eta) \in \Z[a_1,a_2,\ldots, a_m, b_1, b_2, \ldots, b_n]
\]
is a multivariable polynomial, of degree $d$, with $\deg(a_i)=\deg(b_i)=i$, satisfying the definitions below.

\begin{definition}[Intuitive definition]
\label{def:Tp_intuitive}
    For ``nice'' maps (discussed below) $f:M^m\to N^n$ between complex compact manifolds the fundamental class of the closure of 
    \[
\eta(f)=\{
p\in M : \text{the germ of $f$ at $p$ belongs to } \eta\} 
\]
in $H^{2d}(M)$ equals the value of $\Tp(\eta)$, when we substitute $a_i=c_i(TM)$, $b_i=f^*(c_i(TN))=c_i(f^*(TN))$.
\end{definition}

\begin{definition}[More precise definition]
    \label{def:Tp_more_precise}
    The Thom polynomial of $\eta$ is the $\E^o(m)\times \E^o(n)$-equivariant fundamental class of $\overline{\eta} \subset \E(m,n)$. 
\end{definition}

A lot of explanations are in order. Let us start with Definition \ref{def:Tp_more_precise}. The group $\E^o(m)\times \E^o(n)$ and the space $\E(m,n)$ are infinite dimensional, yet via jet approximation (cf. Remark~\ref{rem:jets}) we deal with them as if they were finite dimensional. However, even in finite dimensional settings the concept of equivariant fundamental class is not an obvious one---note for example that $\E(m,n)$ is not compact, only ``equivariantly compact''. A number of precise definitions exist for ``equivariant fundamental class'' \cite{Kaz97,EG98,FR04,MS,AF}, and each one shares key properties with ordinary fundamental classes: 
\begin{itemize}
 \item its cohomological degree is $2d$;
    \item restricted to an invariant set disjoint from $\overline{\eta}$ is 0;
    \item restricted to the smooth part of $\eta$ it is the (equivariant) Euler class of its normal bundle;
    \item behaves ``nicely'' with respect to pullback and push-forward.
\end{itemize}
Accepting the notion of equivariant fundamental class, Definition \ref{def:Tp_more_precise} puts $\Tp(\eta)$ in the ring
\[
H^*_{\E^o(m)\times \E^o(n)} \E(m,n).
\]
The space $\E(m,n)$ is (equivariantly) contractible, so it can be replaces with just a point, and it is a fact that $\E^o(m)$ is homotopy equivalent to $\GL_m(\C)$. Hence we have that $\Tp(\eta)$ is a degree 2d element in
\[
H^*_{\GL_m(\C)\times \GL_n(\C)}(\pt)=
H^*( B\!\GL_m(\C)\times B\!\GL_n(\C))=
\Z[a_1,a_2,\ldots,a_m,b_1,b_2,\ldots,b_n],
\]
the declared habitat of the Thom polynomial.

In the rest of this section we connect Definitions \ref{def:Tp_intuitive} and \ref{def:Tp_more_precise}. This argument is called the {\em degeneracy loci interpretation of Thom polynomials}. The argument is walking through the chain of equalities:
\begin{multline}
\label{eq:Tp_def_chain}
\left[\overline{\eta(f)}\subset M\right]
=
\left[j_f^{-1}\left( \overline{\eta(E)}\right) \subset M\right]
=
j_f^*\left( \left[\overline{\eta(E)}\subset E\right] \right)
\\
=\left[\overline{\eta(E)}\subset E\right]
=k^*_f\left( \left[\overline{\eta}\subset \E(m,n)\right] \right)=
[\overline{\eta}\subset \E(m,n)]|_{a_i=c_i(TM), b_i=c_i(TN)}.
\end{multline}
Here is the explanation: to the map $f:M^m\to N^n$ we associate a fiber bundle $E\to M$ by letting the fiber over $p\in M$ be 
\[
\{\text{germs } (M,p) \to (N,f(p))
\}=
\{\text{germs } (U,p) \to (V,f(p))
\}\cong
\E(m,n),
\]
where $U$ and $V$ are affine neighborhoods of $p$ and $f(p)$. The structure group of this bundle is $\E^o(m)\times \E^o(n)$. Since $\eta\subset \E(m,n)$ is invariant under this group, $\eta$ as a set is defined in every fiber. We call the union of these $\eta$'s $\eta(E)$. Moreover the map $f$ induces an obvious section $j_f$ of the bundle $E\to M$. 

The first equality in \eqref{eq:Tp_def_chain} is the observation that $\eta(f)=j_f^{-1}(\eta(E))$, by definition. The second equality follows from the consistency of fundamental class and pull-back. But let us stop here: this consistency holds, if $j_f$ is transversal to $\eta(E)$ in the right sense---this transversality is exactly the requirement for the map $f$ to be ``nice'' in Definition \ref{def:Tp_intuitive}. The third equality is just a topological fact: the total space of a bundle with contractible fiber is homotopy equivalent to the base. This fact identifies their cohomologies, and any section induces the identity between these cohomologies.
The rest of the equalities follow from the pull-back diagram
\[
\tikz{
\draw (0,0) node {$\eta(E) \subset E$};
\draw[->] (0,-.5) -- (0,-1);
\draw (0,-1.5) node {$M$}; 
\draw (-1,0.5) to[out=-100,in=100] (-1,-2);
\draw (1,0.5) to[out=-80,in=80] (1,-2);
\draw[->] (2,-.8) -- (3,-.8);
\draw (8,0) node {$\eta\subset \E(m,n)$};
\draw[->] (8,-.5) -- (8,-1);
\draw (8,-1.5) node {$point$}; 
\draw (7,0.5) to[out=-100,in=100] (7,-2);
\draw (9,0.5) to[out=-80,in=80] (9,-2);
\draw (5.7,-.8) node {$B_{\E^o(m)\times \E^o(n)}$};
\draw (2.5,-0.6) node {$k_f$};
}
\]
where the specific situation $\eta(E)\subset E$ is pulled back from the universal situation (the Borel construction applied to $\eta\subset \E(m,n)$), via the classifying map $k_f$ that, on cohomology, maps $a_i\mapsto c_i(TM)$, $b_i\mapsto f^*c_i(TN)$.

\bigskip

\noindent{\bf Notation.} Our notation for Thom polynomials will be $\Tp(\ )$, where in the parenthesis we put enough information that specifies the singularity $\eta$.

\subsection{Dependence on quotient variables only}

The reader, no doubt, would like to see examples for Thom polynomials, and is welcome to glimpse at later pages of this paper. However, to avoid complicated, meaninglessly incomplete formulas like the 30-term expression
\begin{equation}\label{eq:bad}
\Tp(A_1,m=10,n=15)=a_1^6-a_1^5b_1-5a_1^4a_2+\ldots -a_6+b_6,
\end{equation}
first let us start with the following fundamental result.

\begin{theorem}[Thom-Damon-Ronga theorem] \label{thm:DamonRonga}
    The Thom polynomial of the contact singularity corresponding to the local algebra $Q$ depends only on the ``quotient variables'' $c_i$ defined by
    \[
    1+c_1t+c_2t^2+\ldots=
    \frac{1+b_1t+
b_2t^2+\ldots}{1+a_1t+a_2t^2+\ldots}.
    \]
    Moreover, when expressed in the quotient variables, the expression only depends on $\ell$, not on $m$ and $n$ individually.
\end{theorem}

For example, statement \eqref{eq:bad} can now be stated concisely as 
\[
\Tp(A_1,m=10,n=15)=c_6,
\]
and even better, the {same} holds for any $n=m+5$, and we just write (cf. \eqref{eq:A1} below)
\[
\Tp(A_1,\ell=5)=c_6.
\]

\subsection{First examples}
The triviality $\Tp(\Sigma^0)=1$ just encodes the fact that a generic map is non-singular at almost all of the points of the source. The first non-trivial general Thom polynomial example is 
\begin{equation}\label{eq:A1}
\Tp(\Sigma^1,\ell)=\Tp(A_1,\ell)=c_{\ell+1}.
\end{equation}

To put this statement in context let us recall the following argument from differential topology: if $c_{\ell+1}(f)=c_{\ell+1}(f^*TN-TM)$ is not $0$ then $f$ must have singularities, it cannot be an immersion. Indeed, if $f$ was an immersion then $f^*TN-TM$ would be a rank $\ell$ {\em bundle}, not just a {\em virtual bundle}---and hence its $\ell
+1$'st Chern class would be 0. The $\Tp(\Sigma^1,\ell)=c_{\ell+1}$ statement makes this ``qualitative'' statement ``quantitative'': $c_{\ell+1}(f)$ is exactly the fundamental class of the locus where $f$ is not an immersion.

Our next example, the Giambelli-Thom-Porteous formula \cite{thomoriginal, porteous}, is also a key formula of Schubert calculus:
\begin{equation*}\label{eq:GTP}
\Tp(\Sigma^r,\ell)=
\det
\begin{pmatrix}
c_{r+\ell} & c_{r+\ell+1} & c_{r+\ell+2} & \ldots  \\
c_{r+\ell-1} & c_{r+\ell} & c_{r+\ell+1} & \ldots \\
 & & & \ddots \\
 & & & \ldots &  c_{r+\ell}  & c_{r+\ell+1}  \\
 & & & \ldots & c_{r+\ell-1} & c_{r+\ell}
\end{pmatrix}_{r\times r},
\end{equation*}
where the name of the right hand side is {\em Schur polynomial} $s_{(r+\ell)^r}$ (cf. Section~\ref{sec:Schur}).

Here are some random other Thom polynomials:
\begin{eqnarray*}
    \Tp(A_2,\ell=0) & = & c_1^2+c_2 \\
    \Tp(A_2,\ell=1) & = & c_2^2+c_1c_3+2c_4 \\
    \Tp(A_3,\ell=0) & = & c_1^3+ 3c_1c_2 +2c_3 \\
    \Tp(I_{22},\ell=1) & = & 2 c_{3} c_{4}-2 c_{2} c_{5}+ c_{1} c_{3}^2-c_{1} c_{2} c_{4},
\end{eqnarray*}
and the interested reader will find other concrete Thom polynomials on the Thom Polynomial Portal \cite{TPP}.


\subsection{Methods of calculating Thom polynomials}
Nobody knows how to calculate the Thom polynomial of a random singularity. Nevertheless, some more-or-less effective methods exist. They are based on one or more of the following geometric tools: resolution (or embedded resolution), partial resolution, degeneration, eg. Gr\"obner degeneration, interpolation, quotient constructions, Hilbert schemes. 

It is also a remarkably successful method to collect {\em structure theorems} on Thom polynomials, and hope that those structure theorems leave very little ambiguity in the coefficient of a sought Thom polynomial. We already met one such structure theorem, Theorem \ref{thm:DamonRonga}.

In the rest of the paper we will meet some computational methods, and some structure theorems. We start with the illustration of the  `interpolation method'.

\subsection{Thom polynomial of $A_4$, $\ell=0$, via interpolation}
\label{sec:interpolation}
Since the codimension of $A_4$ for $\ell=0$ is $4$, we know that 
\begin{equation}\label{eq:Tp_unknown_coeffs}
\Tp(A_4, \ell=0)=Ac_1^4 + B c_1^2c_2 + 
C c_1c_3 + D c_2^2+E c_4
\end{equation}
for some integers $A,B,C,D,E$. The fact is that we can calculate $A,B,C,D,E$ as the unique solution of the system of equations:

\begin{center}
\begin{tabular}{llllll}
     $A(a)^4$ &  $+B(a)(-a^2)$ & $+C(a)(a^3)$ & $+D(-a^2)^2$ & $+E(-a^4)$ & $=0$ \\ 
     $A(2a)^4$ &  $+B(2a)(-2a^2)$ & $+C(2a)(2a^3)$ & $+D(-2a^2)^2$ & $+E(-2a^4)$ & $=0$ \\ 
          $A(3a)^4$ &  $+B(3a)(-3a^2)$ & $+C(3a)(3a^3)$ & $+D(-3a^2)^2$ & $+E(-3a^4)$ & $=0$ \\ 
               $A(4a)^4$ &  $+B(4a)(-4a^2)$ & $+C(4a)(4a^3)$ & $+D(-4a^2)^2$ & $+E(-4a^4)$ & $=24a^4$ 
\end{tabular}
\end{center}
\begin{multline*}
A(a + b)^4 + B(a + b)^2(-a^2 + ab - b^2) + C(a + b)(a^3 - a^2b - ab^2 + b^3) +\\
D(-a^2 + ab - b^2)^2+
E(-a^4 + a^3b + a^2b^2 + a b^3 - b^4) =0.
\end{multline*}
The first four equations are in $\Z[a]$ and the last one is in $\Z[a,b]$. Together they form a system of linear equations for $A,B,C,D,E$ with the unique solution $A=1$, $B=6$, $C=9$, $D=2$, $E=6$, concluding
\begin{equation}\label{eq:A4Tpresult}
\Tp(A_4, \ell=0)=c_1^4 + 6 c_1^2c_2 + 
9 c_1c_3 + 2 c_2^2+6 c_4.
\end{equation}
How did we get the five equations? Let us explain the 4th one. It comes from applying the defining property of the Thom polynomial
\begin{equation}\label{eq:Tp_def}
   \Tp( c(f) )=[\overline{A_4}(f)] 
\end{equation}
to the map $f:\C^4 \to \C^4$
\begin{equation}\label{eq:A4proto}
f:(x,y_1,y_2,y_3)\mapsto (x^5+y_3x^3+y_2x^2+y_1x, y_1, y_2, y_3),
\end{equation}
well, in a sophisticated way. Namely, first notice that the map \eqref{eq:A4proto} has a $U(1)$ symmetry: multiplying the coordinates of the source by $(\alpha, \alpha^4, \alpha^3, \alpha^2)$ multiplies the coordinates of the target by $(\alpha^5, \alpha^4, \alpha^3, \alpha^2)$, for $\alpha\in U(1)$. In effect, we can consider \eqref{eq:Tp_def} is $U(1)$-equivariant cohomology. The right hand side is the class of $0$ in $\C^4$ (with the action $(\alpha, \alpha^4, \alpha^3, \alpha^2)$), which is hence the Euler class $a(4a)(3a)(2a)=24a^4$ of the representation. The left hand side is \eqref{eq:Tp_unknown_coeffs} with $c_i$ being the Taylor coefficients of 
\begin{multline*}
c(f)=\frac{(1+5a)(1+4a)(1+3a)(1+2a)}{(1+a)(1+4a)(1+3a)(1+2a)}=
\frac{1+5a}{1+a}=1+\underbrace{4a}_{c_1}\underbrace{-4a^2}_{c_2}\underbrace{+4a^3}_{c_3}\underbrace{-4a^4}_{c_4}+\ldots.
\end{multline*}
We obtain exactly the non-homogeneous equation listed above. The other, homogeneous, equations are obtained similarly, by applying \eqref{eq:Tp_def} to maps not having $A_4$-points at all: namely to the `prototypes' of the singularities $A_1$, $A_2$, $A_3$, and $I_{22}$ (with their symmetries):
\[
\begin{array}{rl}
  (x)  \mapsto & (x^2) \\
  (x,y) \mapsto & (x^3+xy, y) \\
  (x,y_1,y_2) \mapsto & (x^4+y_2x^2+y_1x, y_1, y_2) \\
  (x,y,u,v) \mapsto & (x^2+uy, y^2+vx,u,v).
\end{array}
\]

The method just described is called the ``interpolation'' method \cite{rrtp,rrA4}. It depends on understanding the symmetries of singularities, and it reduces the calculation of a Thom polynomial to solving a system of linear equations (alternatively, finding the lowest degree generator of the intersection of a bunch of ideals in a polynomials ring).

When this method works, it is surprisingly effective. It has, however, limitations: it can be shown to work for a singularity $\eta$ if all the singularities $\zeta$ with $\codim(\zeta)\leq \codim(\eta)$ satisfy a non-vanishing Euler class condition. This condition holds for all singularities until the appearance of moduli. For example, for $\ell=0$ the interpolation method can calculate the Thom polynomials of exactly the codim $\leq 8$ singularities in Figure \ref{Fig:l0}.

It is important to point out another limitation of the interpolation method: it considers the Thom polynomials for the same algebra but different $\ell$ completely different problems. Hence it is not well positioned to discover structures we will start describing in the next section.

\section{Thom series}

\subsection{Thom polynomials for various $\ell$ values.}
\label{sec:Ts}

Consider the Thom polynomial of the contact singularity $A_3$ for $\ell=0,1,2$, respectively \cite{a3}:
\begin{center}
  \begin{tabular}{p{.1cm}p{1.2cm}p{.9cm}p{.6cm}p{1.1cm}p{1.3cm}p{1cm}p{1cm}p{1cm}l}
$c_1^3$ & $+3c_2c_1$ & $+2c_3$, \\
$c_2^3$ & $+3c_3c_2c_1$ & $+2c_4c_1^2$ & $+c_3^2$ & $+7c_4c_2$ & $+10c_5c_1$ & $+12c_6$,\\
$c_3^3$ & $+3c_4c_3c_2$ & $+2c_5c_2^2$ & $+c_4^2c_1$ & $+7c_5c_3c_1$ & $+10c_6c_2c_1$ & $+12c_7c_1^2$
& $+60c_8c_1$ & $+26c_7c_2$ & $+\ldots$.
  \end{tabular}
\end{center}

The pattern becomes clearer if we imagine every monomial to be a product of {\em exactly three} $c_i$ variables---we call this ``three'' the {\em width} of the monomials. For this, we observe that there are no monomilas whose width is $\geq 4$, and the monomial whose width is strictly less than $3$ we formally multiply by the appropriate $c_0$-power ($c_0=1$). 

Then each line is obtained from the line below it by the ``lowering'' $\flat$ operation, where $(c_ic_jc_k)^{\flat}=c_{i-1}c_{j-1}c_{k-1}$. Observe that, for example, $(c_3^2)^{\flat}=(c_3c_3c_0)^{\flat}=c_2c_2c_{-1}=0$, and indeed, 
\[
 \Tp(A_3,\ell)^{\flat}=\Tp(A_3,\ell-1)
\]
holds for $\ell>1$.

It is tempting to define the raising operator $(c_ic_jc_k)^{\#}=c_{i+1}c_{j+1}c_{k+1}$, and expect that $\Tp(A_3,\ell)^{\#}=\Tp(A_3,\ell+1)$, but---as the examples above show---this idea determines only {\em some} of the monomials of $\Tp(A_3,\ell+1)$, exactly those whose $\flat$-image is not 0. 

The right way of encoding the structure we just illustrated is to consider the limit object $\Tp(A_3,\ell=\infty)$, called Thom series for $A_3$, (with shifted indices) as follows:
\[
\Ts(A_3)=
d_0^3+3d_1d_0d_{-1}+2d_2d_{-1}^2+
d_1d_{-2}^2+7d_2d_0d_{-2}+10 d_{3}d_{-1}d_{-2}+\ldots
\]
Then we have
\[
\Tp(A_3,\ell)=\Ts(A_3)|_{d_i=c_{i+\ell+1}}.
\]
The point is that $\Tp(A_3,\ell)$ for all $\ell$ is determined by just {\em one} formal power series $\Ts(A_3)\in \Z[\![\ldots,d_{-2},d_{-1},d_0,d_1,d_2,\ldots]\!]$.

\smallskip

The described $A_3$ example, namely, the existence of one Thom series 
\[
\Ts(\eta)\in \Z[\![\ldots,d_{-2},d_{-1},d_0,d_1,d_2,\ldots]\!]
\]
from which all $\Tp(\eta,\ell)$'s are obtained by the $c_i=d_{i+\ell+1}$ substitution holds for a wide class of contact singularities and Thom-Boardman singularities (and is believed to hold for all)~\cite{dstab}.

The width and the degree of the Thom series are easily identifiable invariants of $\eta$. For example, if $\eta$ is a contact singularity $\eta_Q$, then the width is $\dim(Q)-1$.

\subsection{Thom series in Schur expansions}\label{sec:Schur}
  It is often worth expressing Thom polynomials not in Chern monomial basis, but in the basis of Schur polynomials: 
  \[
  s_{\lambda}=\det(c_{\lambda_i+j-i})_{i,j} \qquad\qquad (\text{where }\lambda \text{ is a partition}).
  \]
  For example 
\begin{equation}\label{eq:A2Schur}
  \begin{array}{ll}
       \Tp(A_2,\ell=0)& =s_{11}+2s_{2} \\
       \Tp(A_2,\ell=1)& =s_{22}+2s_{31}+4s_{4}, \\
       \Tp(A_2,\ell=2)&=s_{33}+2s_{42}+4s_{51}+8s_{6}.
  \end{array}
\end{equation}
For such expressions the Thom series phenomenon takes this form: there exists a formal expression 
  \begin{equation} \label{eq:SchurTs}
  \Ts(\eta)=\sum_\lambda p_\lambda r_\lambda
  \end{equation}
  where $\lambda$'s are weakly decreasing {\em integer} sequences of a fixed length (the ``width'' above), $p_\lambda\in \Z$, and for any specific $\ell$ we have 
  \[
  \Tp(\eta,\ell)=\Ts(\eta)^{\#(l+1)}:= \sum_\lambda p_\lambda 
  r_{\lambda_1+\ell+1, \lambda_2+\ell+1,\ldots}.
  \]
The Thom series of $A_2$ in this Schur form is given by:
  \[
  \Ts(A_2)=r_{00}+2r_{1,-1}+4r_{2,-2}+8r_{3,-3}+\ldots.
  \]
As the reader may have noticed, we use the same notation for the Thom series regardless of whether it is expressed as a formal power series in $d_i$ variables or in the Schur form \eqref{eq:SchurTs}.

\bigskip

In the next three sections we will describe three different ways of encoding the whole Thom series by just ``finite information.'' Yet, even that finite information is very challenging to study. For example, the Thom series of $A_7$ in unknown.

\subsection{Finite encoding of Thom series \`a la \cite{BSz}}
\label{sec:finite1}

A remarkable new method for presenting Thom series was found in \cite{BSz}, with a clarifying re-interpretation in \cite{kaza:noas}. This approach constructs partial resolutions of singularity loci using the ``test curve method" or through the ``space of non-associative algebras." Equivariant localization and its residue formalism are applied to the partial resolution to obtain formulas for Thom series. For example 
\begin{multline*}
\Ts(A_3)=
-
\mathop{\Res}_{z_1=\infty}
\mathop{\Res}_{z_2=\infty}
\mathop{\Res}_{z_3=\infty}
\left(
\frac{q(z_1,z_2,z_3) }{(2z_1-z_2)(2z_1-z_3)(z_1+z_2-z_3)} 
\right. \times
\\
 \prod_{1\leq i<j\leq 3}(z_i-z_j)
\cdot 
\left.
\left(\prod_{i=1}^3 \sum_{j=-\infty}^{\infty} \frac{d_{j}}{z_i^j}\right)
\frac{dz_3}{z_3}
\frac{dz_2}{z_2}
\frac{dz_1}{z_1}
\right),
\end{multline*}
where $q(z_1,z_2,z_3)=1$. 

For other algebras $Q$ the nature of the formula is similar: The number of the auxiliary $z_i$ variables is $\dim(Q)-1$. The factors in the second line generalize the obvious way. The denominator of the fraction in the first line is also a more-or-less obvious invariant of $Q$, it is always a product of linear factors. The only non-obvious ingredient in general is the $q$-polynomial in the numerator. (The trivial value $q=1$ for $Q=A_3$ is misleading.)

It is hence the $q$-polynomial that is the ``finite'' encoding of the whole Thom series. The polynomial $q$ has a geometric meaning, it is itself an equivariant fundamental class of a variety in a finite dimensional vector space. The vector space is the space of (not necessarily associative) algebras with some numerical characteristics. The subvariety whose fundamental class is $q$ is the collection of points corresponding to associative algebras isomorphic to~$Q$.

In effect, the calculation of the Thom series of $Q$ is reduced to understanding how algebras isomorphic to $Q$ sit in the space of algebras. While this reduction is satisfying, it is discouraging that---in general---we do not know algebras very well. For a general $Q$ we do not know the $q$ polynomial. 

In small examples, however, $q$ can be computed. For example for $Q=A_5=\C[x]/(x^6)$ we have
\[
q=(z_5-2z_1-z_2)(2z_1^2+3z_1z_2-z_1z_5+zz_2z_3-z_2z_4-z_2z_5-z_3z_4+z_4z_5),
\]
but the $q$-polynomial of, say, $A_{10}$ seems way beyond the scope of our present understanding of algebra/geometry.

\begin{remark}
   The iterated residue form above may look complicated at the first sight. In fact, it is very explicit. One takes the appropriate expansion $\sum_{ijk} a_{uvw} z_1^uz_2^vz_3^w$ of the rational function. The effect of multiplying with $\prod \sum d_{j}/z_i^j$ and taking residues is that the monomial $z_1^uz_2^vz_3^w$ turns to $d_ud_vd_w$. Hence, the formula above almost identifies $\Ts(A_3)$ with the expansion of a rational function. The simple but subtle difference is that terms with permuted $u,v,w$-exponents give the same $d$-monomial. 
\end{remark}

\begin{remark}
     When setting up the ``space of algebras''  one has the freedom of {\em choosing} a filtration on the algebra $Q$. For any such filtration a residue form of the Thom series is obtained (with the non-obvious ingredient $q$). It is a fact that for some filtration choices the $q$-polynomial is much easier to find than for others, see examples in \cite{kaza:noas, TPP}. We cannot illustrate this phenomenon with $A_3$, which has only one filtration. 
\end{remark}

\begin{remark}
    We have not provided a {\em precise} definition of the $q$-polynomial or its geometric origin. Therefore, the reader must take on faith that its geometry is governed by a Borel group action. The non-reductive nature of Borel groups presents a significant challenge in advancing this approach further (cf.~Section~\ref{sec:finite3}). 
\end{remark}

\subsection{Finite encoding of Thom series a la \cite{FRannals}}
\label{sec:finite2}

In \cite{FRannals}
a partial geometric resolution (pioneered by J. Damon \cite{damon}) and singular equivariant localization methods  are used to reduce the Thom series $\Tp_Q$ to finitely many rational functions. 

Let us illustrate the nature of such formulas by an example. Let $A=\{a_1,\ldots,a_m\}$, $B=\{b_1,\ldots,b_n\}$, and define $\RE(X,Y)=\prod_{y\in Y, x\in X} (y-x)$. We claim
\[
\Tp(A_2,m \to n)=
\sum_{i=1}^n
\frac{\RE(W_1,B)}{e_1 \cdot \RE(W_1,A_1)}
+
\sum_{1\leq i < j\leq n}
\frac{\RE(W_2,B)}{e_2 \cdot \RE(W_2,A_2)},
\]
where $
W_1=\{\alpha_i,2\alpha_i\},
A_1=A-\{a_i\}$ and $
W_2=\{\alpha_i,\alpha_j\},
A_2=A-\{a_i,a_j\}$, and, most notably
\[
e_1=1, \qquad\qquad e_2=\frac{1}{3}(a_1-2a_2)(a_2-2a_1).
\]
Some explanations are in order. The two summations correspond to the two possible {\em monomial ideal types} of length 3, namely: $(x^3) \triangleleft \C[x]$ and $(x^2,xy,y^2) \triangleleft \C[x,y]$. Some ingredients of the formula are obvious: $W_i$ are torus weights on the quotient by the monomial ideal, $A_i$ just record the $a$ variables {\em not used} in the monomial ideal. The only non-obvious ingredient, coming from geometry, are the classes~$e_i$.

The example we gave generalizes to any algebra $Q$. If $\dim(Q)=\mu+1$ the formula breaks into terms corresponding to the length $\mu+1$ monomial ideal types. In these terms all but the rational functions $e_i$ are obvious. In summary: the Thom polynomials $\Tp(Q,m\to n)$ for all $m,n$, equivalently the Thom series $\Ts$, is determined by a finite sequence of rational functions~$e_i$.

On the one hand, the $e_i$ classes have geometric meaning, on the other hand, they satisfy intriguing algebraic relations \cite{FRannals}. The interplay of the two worlds makes the $e_i$ (and hence the Thom polynomials $\Tp(Q,\ell)$ for all $\ell$) calculable for algebras of small dimension.

An interesting consequence  is the technique of reverse engineering: if we know $\Tp(Q)$ for {\em some} $m,n$ explicitly (say, by interpolation, see Section \ref{sec:interpolation}), we can use that information to find all $e_i$, and then we have the Thom polynomial for all $m,n$.

\subsection{Finite encoding of Thom series a la \cite{GBnr}}
\label{sec:finite3}
An improvement of the method of Section \ref{sec:finite1} is presented in \cite{GBnr}. As we pointed out, the geometry behind that approach is a non-reductive group action. Recent advances on understanding non-reductive actions and quotients lead G. B\'eczi to modify the resolution based on ``test-curves'' to another one which is {\em toric} in nature. The price to pay is that the resolution is now an iterated one, and the structure of the iteration is arranged in a tree. The key ingredient of the obtained formula is a sum whose terms correspond to the leaves of the tree. Each such contribution is of very simple structure---one only needs to work with products of linear factors as opposed to general polynomials. For more details we refer the reader to the paper \cite{GBnr}, here we only show one example:
\begin{multline*}
    \Ts(A_4)=
\mathop{\Res}_{z_1=\infty}
\mathop{\Res}_{z_2=\infty}
\mathop{\Res}_{z_3=\infty}
\mathop{\Res}_{z_4=\infty}
\Biggl(
{G(z_1,z_2,z_3,z_4)}\times 
 \\
\left.
\prod_{1\leq i < j \leq 4} (z_i-z_j) \cdot
\left(\prod_{i=1}^4 \sum_{j=-\infty}^{\infty} \frac{d_{j}}{z_i^j}\right)
\frac{dz_4}{z_4}
\frac{dz_3}{z_3}
\frac{dz_2}{z_2}
\frac{dz_1}{z_1}
\right),
\end{multline*}
where
\[
G=\frac{
\frac{1}{(2z_1-z_2)(z_1+z_2-z_3)} - \frac{1}{(z_1+z_3-z_4)}}{(z_2+z_3-z_1-z_4)(z_1+z_2-z_4)(2z_1-z_3)(2z_1-z_4)}.
\]

\section{Real singularities}

In this section we consider maps $f: M^m \to N^n$ between real manifolds, and correspondingly real singularities $\eta \subset \{(\R^m,0)\to (\R^n,0)\}$. We are interested in Thom polynomials expressing the fundamental class of the closure of $\eta$-points of $f$, in terms of the characteristic classes of $f$. 

The first new phenomenon is that real singularity loci do not necessarily carry a fundamental class. 
One of the geometric reasons is that $\eta$-points of a map are just semi-algebraic sets, like a semicircle: the topological closure is not necessarily a cycle.

\subsection{Cohomology with $\Z_2$  coefficient, Stiefel-Whitney classes}

If we aim at Thom polynomials in $\Z_2$-coefficient cohomology, then the mentioned difficulty is essentially the only one. This is the content of the Borel-Haefliger theorem \cite{borel-haefliger}. Namely, if 
\[
\eta^{\R} \subset \{(\R^m,0)\to (\R^n,0)\}
\qquad\text{ and } \qquad  
\eta^{\C} \subset \{(\C^m,0)\to (\C^n,0)\}
\]
are the real and complex `forms' of each other, then the $\Z_2$ coefficient Thom polynomial of $\eta^{\R}$ is obtained from the Thom polynomial of $\eta^{\C}$ by replacing Chern classes $c_i$ with the corresponding Stiefel-Whitney classes $w_i$, and reducing the coefficients mod 2. For example, since the Thom polynomial 
\[
\Tp(\Sigma^3, \ell=0)=\det\!
\begin{pmatrix} 
c_3 & c_4 & c_5 \\
c_2 & c_3 & c_4 \\
c_1 & c_2 & c_3
\end{pmatrix}
=
c_3^3+c_4^2c_1+c_5c_2^2-2c_4c_3c_2-c_1c_3c_5,
\]
we have that the $\Z_2$-coefficient Thom polynomial of the $\Sigma^3$ real singularity is 
\[w_3^3+w_4^2w_1+w_5w_2^2+w_1w_3w_5.\]

\subsection{Cohomology with $\Z$ or $\Q$ coefficients, Pontryagin classes \cite{ronga,FRint,CsSzT}}
We meet a richer structures if we want to find Thom polynomials of real singularities with $\Q$ or even $\Z$ coefficients. Here is an example \cite{FRint}: 
\[
\Tp(\eta_1 - 2\eta_2,\ell=0)=
\left(
3p_1^2+9p_2 
\right)+
\left(
v_1^2v_2v_4 + 
v_1v_2v_5 + v_1v_3v_4+v_3v_5
\right).
\]
The right hand side is a characteristic class of $f$ (that is, of $f^*TN-TM$) as before, but now the classes are not Chern classes. The $p_i$ denote Pontryagin classes and the $v$-monomials $v_iv_jv_k$ are the unique integer lifts of the Stiefel-Whitney monomials $w_iw_jw_k$. These $v$-monomials are 2-torsion elements in $\Z$-coefficient cohomology. The $\Q$-coefficient Thom polynomial is hence $3p_1^2+9p_2$.

What needs more explanation is the left hand side: we stated the Thom polynomial of what? The answer is that our cycle is the (closure of the) singularity locus for $\eta_1=\R[x,y]/ (x^2+y^3,xy^2)$ (with a particular orientation, not described here) {\em minus} twice the singularity locus  for $\eta_2=\R[x,y]/(x^2+y^2,x^4)$ (with a particular orientation, not described here). 

Several non-trivial statements are encoded in the last sentence. First: these two singularity loci have 8 real codimension and are both co-orientable. If a singularity locus in not co-orientable, then even if it had no boundary (at more complicated singularities) it would not carry an integer coefficient cohomology class in an oriented manifold---this is another of the geometric difficulties we alluded to above. 
\begin{figure}
\[
\begin{tikzpicture}[
    middlearrow/.style 2 args={
        decoration={             
            markings, 
            mark=at position 0.5 with {\arrow[xshift=3.333pt]{triangle 45}, \node[#1] {#2};}
        },
        postaction={decorate}
    },
]
  \node (A) at (0,0)   {$\bullet$};
  \node (B) at (2.7,2)   {};
  \node (C) at (2.7,1.4)   {};
  \node (D) at (2.7,.6) {};
\draw[solid,fill] (A) circle (.15);
\draw[middlearrow,thick]   (A) to[out=60,in=180] (B);
\draw[middlearrow,thick]   (A) to[out=50,in=180] (C);
\draw[blue] (3.2,2.1) to[out=0,in=180] (3.5,1.7);
\draw[blue] (3.2,1.3) to[out=0,in=180] (3.5,1.7);
\node[blue] at (5.7,1.7) {$\eta_1=\R[x,y]/(x^2+y^3,xy^2)$};
\draw[middlearrow,thick]   (A) to[out=30,in=170] (D);
\node[blue] at (5.3,.6) {$\eta_2=\R[x,y]/(x^2+y^2,x^4)$};
\node[blue] at (-2,1) {$\R[x,y]/(x^3+y^3,y^4)$};
\draw[blue,->] (-2,.8) to[out=-90,in=180] (-0.5,0);
\node at (2.85,2) {$\cdots$};
\node at (2.85,1.4) {$\cdots$};
\node at (2.85,.6) {$\cdots$};
\end{tikzpicture}
\]
\caption{The picture of the 8-codimensional singularities $\Sigma_1$ and $\Sigma_2$ in a normal slice to a 9-codimensional singularity.}
\label{fig:VassCompGeo}
\end{figure}
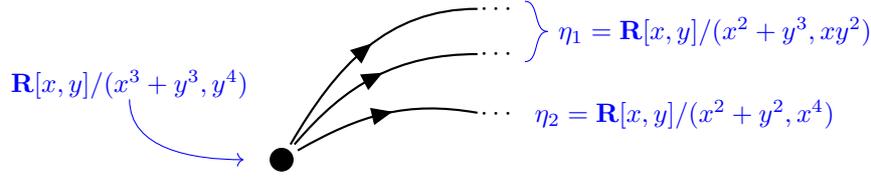
Second: the local picture of $\eta_1$- and $\eta_2$-loci near the 9-codimensional singularity $\R[x,y]/(x^3+y^3,y^4)$ is illustrated in Figure \ref{fig:VassCompGeo}. Accordingly, only (the closure of) $\eta_1-2\eta_2$ (or its integer multiples) forms a cycle around this point, no other linear combination of $\eta_1$ and $\eta_2$. Third: there are no other codimension 9 singularities besides $\R[x,y]/(x^3+y^3,y^4)$ in the boundary that would make $\eta_1-2\eta_2$ not a cycle. 

The geometric information needed to find similar linear combinations of real singularities is organized into the so-called Vassiliev algebraic complex \cite{Kaz97}, a cousin of the famous Vassiliev complex in knot theory.

\section{Generalized Thom polynomial concepts}

\subsection{The landscape of charactersitic classes}
\label{sec:landscape}

According to Definition \ref{def:Tp_more_precise}, the Thom polynomial of a singularity $\eta$ is an (equivariant) {\em cohomological fundamental class} of the subvariety $\overline{\eta}$ in a vector space. 

Some notable generalizations of the  concept {\em cohomological fundamental class} are arranged in the table
\begin{equation}\label{table:generalized}
    \begin{tabular}{c|c|l}
         & fundamental class & $\hbar$-deformed fundamental class\\
         \hline
      $H^*$   &  $\exists !$ & Chern-Schwartz-MacPherson class\\
      $K$   &$ \exists \exists \exists$ &  motivic Chern class\\
      $\Ell$ & $\nexists$ & elliptic class (with K\"ahler variables).
    \end{tabular}
\end{equation}
That is, we consider the three generalized cohomology theories whose associated formal group laws are 1-dimensional algebraic groups: cohomology, K theory, and elliptic cohomology. 

It turns out that in K theory there are three non-equivalent ``fundamental class'' concepts: (i) the class of the structure sheaf, (ii) the class of the structure sheaf of a resolution, pushed forward, and (iii) well, a third one we will meet soon. It is a theorem that in elliptic cohomology there is no fundamental  class concept that does not depend on choices. 

Remarkably, there is a one-parameter ($\hbar$) deformation of the fundamental class concept, in all three theories. The cohomological (called Chern-Schwartz-MacPherson, CSM) class, and the K theoretic (called motivic Chern, MC) class were defined classically, as classes satisfying an additivity (motivic) property, and also being consistent with push-forward maps. 

At certain asymptotics the CSM class recovers the cohomological fundamental class. The $\hbar=1$ substitution of the MC class is the ``third'' version of K theoretic fundamental class. The three K theoretic fundamental classes of a singular set $\overline{\eta}$ coincide if $\overline{\eta}$ has only ``mild'' singularities. At the $K\to H^*$ reduction all three reduce to the cohomological fundamental class. 

Motivated by the works of Okounkov, Maulik, Aganagic, which relate characteristic classes to Geometric Representation Theory, as well as definitions of Borisov-Libgober, the ``right" concept of the $\hbar$-deformed elliptic fundamental class has been defined in certain contexts \cite{RWell, KRWell}. Unfortunately, these contexts do not yet encompass Thom polynomials. Furthermore, even in these contexts, the $\hbar$-deformed elliptic class inherently depends on a new set of parameters known as Kähler (or dynamical) parameters.

\begin{remark}
    It is tempting to think that an elliptic fundamental class could be defined by specializing the new parameters and $\hbar$. However, it turns out that at that specialization the $\hbar$-deformed elliptic classes have poles. This fact is another incarnation of the result mentioned above on the non-existence of elliptic fundamental class.
\end{remark}

Suppose we know the $\hbar$-deformed classes for the singular set $\overline{\eta}$ in a smooth ambient space. What new information do these classes provide? They determine the Euler characteristics (for CSM class), the $\chi_y$-genus (for MC class), and the elliptic genus (for elliptic class) of various sections of $\overline{\eta}$.

\smallskip

The question relevant to us is: do the generalizations of Table \eqref{table:generalized} exist for Thom polynomials? If so, can they be computed and applied in geometry? We summarize the currently known partial answers in the next few sections.

\subsection{K theory Thom polynomials \cite{RSz}}
\label{sec:KTp}

The simplest K theory Thom polynomials are objects in ``K theory Schubert Calculus". Namely, the K theory Thom polynomials---denoted by $\KTp$---are known for the Thom-Boardman singularities of order 1:
\[
\KTp(\Sigma^r,\ell)=g_{\underbrace{r+\ell,r+\ell,\ldots,r+\ell}_{r}}.
\]
Here $g_{\lambda}$ is a Grothendieck polynomial, the K theoretic analogue of the Schur polynomial.


Below are the K theory Thom polynomials for the $A_2=\C[t]/(t^3)$ singularity for $\ell=0,1,2$ respectively:
\begin{center}
    \begin{tabular}{p{0.3cm}p{.62cm}p{.62cm}p{.62cm}p{.72cm}p{.72cm}p{.72cm}p{.72cm}p{.72cm}p{.72cm}p{.72cm}p{.72cm}p{.72cm}}
    $(g_{11}$ & $+2g_{2})$ & & & $-(2g_{21}$ & $+g_{3})$ & &  & $+(g_{31}),$ &  &  &  \\
    $(g_{22}$ & $+2g_{31}$ & $+4g_4)$ & & $-(2g_{32}$ & $+5g_{41}$ & $+4g_5)$ & & $+(g_{42}$ & $+4g_{51}$ & $+g_6)$ & & $-(g_{61}),$ \\
    $(g_{33}$ & $+2g_{42}$ & $+4g_{51}$ & $+8g_6)$& $-(2g_{43}$ & $+5g_{52}$ & $+12g_{61}$ & $+12g_{7})$& $+(g_{53}$ & $+4g_{62}$ & $+13g_{71}$ & $+\ldots)$ & $-\ldots$.
    \end{tabular}
\end{center}
The structure is worth studying. First, each  partition that occurs as a $g$-subscript has exactly two parts $\lambda=(ij)$ (note that eg. $g_{4}=g_{40}$).

The first parenthesis in each row recovers the corresponding cohomology Thom polynomial, if we formally replace Grothendieck polynomials $g_\lambda$ with Schur polynomials $s_\lambda$, cf. \eqref{eq:A2Schur}. This part displays the ``Thom series'' stabilization of Section \ref{sec:Ts}, and could be encoded by the infinite expression 
\[
r_{00}+2r_{1,-1}+4r_{2,-2}+\ldots,
\]
via $r_{ij}=g_{i+\ell+1,j+\ell+1}$, cf. Section~\ref{sec:Schur}.

However, this stabilization seemingly breaks down at the higher order terms that contain the new geometric information (compared to cohomology) supported on non-top dimensional strata of $A_2 \subset \E(m,m+\ell)$.

The stabilization, however, can be salvaged---because of the fascinatingly combinatorics of Grothendieck polynomials. To illustrate this, let us apply the lowering operation $\flat$ to the second line above. We obtain:
\[
(g_{11}+2g_{20}+4g_{3,-1})
-
(2g_{21}+5g_{30}+4g_{4,-1})
+
(g_{31}+4g_{40}+g_{5,-1})
-
(g_{50}).
\]
As mentioned, $g_{i0}=g_i$. However, while for Schur polynomials $s_{i,<0}=0$, for Grothendieck polynomials we have $g_{i,<0}=g_{i,0}=g_i$. Therefore, the last displayed expression further equals 
\[
(g_{11}+2g_{20}+4g_{3}) - (2g_{21}+5g_{3}+4g_4)
+
(g_{31}+4g_4+g_5)-(g_5),
\]
which, after cancelling terms in different parentheses, recovers the top line in the above table: $\KTp(A_2,\ell=0)$. What we learned is that 
\[
\KTp(A_2,\ell=1)^\flat = \KTp(A_2,\ell=0)
\]
holds, even though the left hand side is a ``non-economical'' way of naming the right hand side.

The described phenomenon holds for all $\ell$. We obtain that there is a K theory Thom series for $A_2$
\begin{multline*}
\KTs(A_2)=
(r_{00}+2r_{1,-1}+\ldots)-
(2r_{10}+5r_{2,-1}+12r_{3,-2}+28r_{4,-3}+\ldots)+\\
(r_{20}+4r_{3,-1}+13r_{4,-2}+38r_{5,-3}+\ldots)- \ldots
\end{multline*}
that determines $\KTp(A_2,\ell)$ for all $\ell$ via $r_{ij}=g_{i+\ell+1,j+\ell+1}$. Moreover, the ``finite information'' result of Section \ref{sec:finite1} also extends to K theory: the coefficients of $\KTs(A_2)$ are the coefficients of the appropriate Laurent expansion of 
\[
\frac{1-2x_2+x_1^2}{x_2-2x_1+x_1^2}.
\]

\begin{remark}
    As we mentioned in Section \ref{sec:landscape} there are three different notions for K theory fundamental classes, and hence, for K theory Thom polynomials. All $\KTp$'s presented in this section are the ``structure sheaves of a resolution pushed forward'' version. 
\end{remark}

\subsection{Segre-Schwartz-MacPherson (SSM) Thom polynomials} \label{sec:ssmTp}

Recall that the concept of fundamental class $[\Sigma\subset M]$ of a possibly singular subvariety $\Sigma$ in $M$ is consistent both with respect to pull-back 
\[
\text{for } 
\qquad
\begin{tikzpicture}[baseline=-10]
\draw node at (0,0) {$f:M \to N$};
\draw node at (.7,-.3) {$\cup$};
\draw node at (.7,-.65) {$\Sigma$};
\end{tikzpicture}
\qquad
\text{ we have }
\qquad
[f^{-1}(\Sigma) \subset M] = f^*([\Sigma \subset N]),
\]
and with respect to push-forward
\[
\text{for } 
\qquad
\begin{tikzpicture}[baseline=-10]
\draw node at (0,0) {$f:M \to N$};
\draw node at (-.2,-.3) {$\cup$};
\draw node at (-.2,-.65) {$\Sigma$};
\end{tikzpicture}
\qquad
\text{ we have }
\qquad
[f(\Sigma) \subset N] = f_*([\Sigma \subset M]),
\]
under the appropriate assumptions on $f$.
The CSM (Chern-Schwartz-MacPherson) class is a 1-parameter deformation of the fundamental class
\[
\csm(\Sigma\subset M)=[\Sigma\subset M] + \text{higher order terms } \in H^*(M)
\]
such that: 
\begin{itemize}
\item[(i)] It is consistent with push-forward (in a sense that involves the Euler-characteristic of fibers of the map).
\item[(ii)] Its SSM (Segre-Schwartz-MacPherson) version
\[ 
\ssm(\Sigma \subset M)
=\frac{\csm(\Sigma\subset M)}{c(TM)}
\in H^{**}(M)
\]
is consistent with pull-back (with appropriate transversality assumption on the map).
\item[(iii)] CSM classes satisfy normalization:
    if $i:\Sigma\subset M$ is a smooth closed subvariety then $\csm(\Sigma \subset M)=i_*(c(T\Sigma))$.
\item[(iv)] The CSM class is defined for all constructible functions on the ambient space, where the 
$\csm(\Sigma\subset M)$ concept corresponds to the indicator function of $\Sigma$.  In this generality CSM classes (hence SSM classes too) are additive: 
\[
\csm(\alpha f + \beta g)=\alpha \csm(f) + \beta \csm(g).
\]
\end{itemize}
In fact, a precise definition of the CSM class should start with property (iv), and in this language CSM is a natural transformation between the functor of constructible functions and  (Borel-Moore) homology. Our version is its Poincar\'e dual in the ambient space.

The listed four properties (i)--(iv) make CSM classes powerful tools of enumerative geometry. For example, from the higher order terms of $\csm(\Sigma\subset M)$ one can calculate the Euler characteristics of various sections of $\Sigma$, see \cite{aluffi}.

The CSM-, or more precisely, the SSM-deformation of Thom polynomials are due to Ohmoto \cite{toru:SMTP}. Various terms of SSM Thom polynomials can be calculated via geometry \cite{toru:SMTP, Nekarda} or an analogue \cite{RV, FRcsm} of the interpolation property of Maulik-Okounkov stable envelopes \cite{MO}. For example, we have 
\begin{multline*}
    \ssmTp(A_1,\ell=0)=
c_{1} 
+ \left(-2 c_{1}^2-c_{2} \right) 
+ \left(3 c_{1}^3+ c_{3}+ 3 c_{1} c_{2}\right)\\
+ \left(-4 c_{1}^4-6 c_{1}^2 c_{2}-7 c_{1} c_{3}+ 3 c_{2}^2-c_{4}\right)+\ldots,
\end{multline*}
or the same in Schur basis
\begin{multline} \label{eq:SSM_example}
    \ssmTp(A_1,\ell=0)=
s_{1} - \left(2 s_{11}+3 s_{2}\right) 
+ \left(3 s_{111}+ 9 s_{21}+ 7 s_{3}\right) \\ - 
\left(4 s_{1111}+18 s_{211}+11 s_{22}+28 s_{31}+15 s_{4}\right)+ \ldots.
\end{multline}

In \cite{FMrealcsm} the construction of the {\em real} version of SSM Thom polynomials is announced.

\subsection{Motivic Segre (MS) Thom polynomials} \label{sec:msTp}

In cohomology the {\em fundamental class} concept $[\Sigma \subset M]\in H^*(M)$ for a singular subvariety $\Sigma\subset M$ is consistent with push-forward and pull-back. As mentioned, in K theory there are two versions for the concept {\em fundamental class}: the class of the structure sheaf, and the class of a structure sheaf of a resolution pushed forward. Roughly speaking, the first one is consistent with pull-back, and the second one is consistent with push-forward. Neither are consistent with respect to the other operation. The two concepts coincide if the singularity is mild (eg. rational).

Hence it is remarkable that an $\hbar$-deformed fundamental class concept exist that is consistent with both push-forward and pull-back (in a certain sense) \cite{BSY, toru:SMTP, FRWmc}. This characteristic class is called the {\em motivic Chern class}, often denoted by $\MC$, $\MC_y$, or $\MC_{\h}$. The subscript is the deformation parameter, which is $y$ classically (cf. Hirzebruch's $\chi_y$ genus), but we prefer $\h$ showing its place in the landscape described in Section \ref{sec:landscape}. 

At $\h=0$ the motivic Chern class is the ``third'' K theoretic fundamental class version we alluded to in Section \ref{sec:landscape}. In general it is different from the two mentioned above \cite{LFmc}, but for $\Sigma$ with mild singularities it coincides with both. 

Here is an informative summary of some properties of motivic Chern classes (cf. the properties of cohomological CSM classes of Section \ref{sec:ssmTp}): 
\begin{itemize}
    \item[(i)] It is consistent with push-forward (in a sense that involves the $\chi_y$-genus of the fibers of the map).
    \item[(ii)] Its MS (motivic Segre) version
    \[
    \MS(\Sigma\subset M)=\frac{\MC(\Sigma\subset M)}{c^K(TM)}\in K(M)[[\h]]
    \]
    is consistent with pull-back (with appropriate transversality assumption of the map). Here $c^K$ stands for total Chern class of a bundle in K theory. For example $c^K(L)=1+\h L^*$ for a line bundle $L$.
    \item[(iii)] The MC classes satisfy normalization:
    if $i:\Sigma\subset M$ is a smooth closed subvariety then $\MC(\Sigma \subset M)=i_*(c^K(T\Sigma))$.
\item[(iv)] The MC class is defined for all morphisms $f:A\to M$ of algebraic varieties, where the 
$\MC(\Sigma\subset M)$ concept corresponds to the inclusion map $\Sigma \subset M$.  In this generality MC classes (hence MS classes too) are motivic: 
\[
\MC(\alpha f + \beta g)=\alpha \MC(f) + \beta \MC(g).
\]
\end{itemize}
In fact, a precise definition should start with property (iv) in the language of mc being a natural transformation between the functors $\Var(-)$ and $K(-)[\h]$. Here $\Var(M)$ is generated
by classes of maps $A \to M$, where $A$ can be singular, and $f$ is not necessarily proper, modulo
additivity relations \cite{BSY}. 

The listed four properties (i)--(iv) make MC classes powerful tools of enumerative geometry, even more so than CSM classes:
The CSM class can be recovered from the MC class by a simple algebraic procedure see eg. \cite[Rem 2.4]{FRWmc}.

\smallskip

Despite their importance, hardly anything is known about {\em motivic Segre Thom polynomials}. Some initial results are summarized in \cite{LFRRmsThom}. A sample result is the motivic Segre Thom polynomial of the {\em closure} of the $A_2$ singularity, namely:
\begin{multline}
\label{MSTP_example}
    \mSTp(\overline{A_2},\ell=0)=
\frac{1}{(1+\h)^2}\Big( g_{11} +  2g_{2} \Big)
\\
-
\frac{1}{(1+\h)^3}
\Big(
 (-2\h)g_{111}+ (2-5\h)g_{21} + (1-5\h)g_{3}
 \Big)\\
+
\frac{1}{(1+\h)^4}
\Big(
 (3\h^2)g_{1111}+ (-7\h+8\h^2)g_{211} + (1-13\h+11\h^2)g_{31}+\ldots\Big) + \ldots.
\end{multline}
(The first coefficient that is unknown is that of $g_{22}$.)
Notice that the $\h=0$ substitution recovers the 
K theoretic Thom polynomial 
\[
\KTp(A_2,\ell=0)=
(g_{11} +2g_{2}) -(2g_{21} +g_{3})+(g_{31})
\]
from Section \ref{sec:KTp}.

\begin{remark}
    The $\h=-1/(s-1)$ substitution in the ms-Tp above gets rid of denominators:
    \begin{multline}\label{eq:MSTP_another}
    \frac{1}{s^2}\mSTp(\overline{A_2},\ell=0)=
    \Big(g_{11}+2g_2\Big)  \\ +\Big((2+3s)g_{111} + (5+12s) g_{21} + (5+12s)g_3\Big)+\ldots,
    \end{multline}
    displaying further positivity (and other) properties of the coefficients. 
\end{remark}

\section{Positivity}

In our summary of various versions of Thom polynomials we focused on their stability structure theorems---such as Theorem \ref{thm:DamonRonga} and Section \ref{sec:Ts}. 
Another intriguing bouquet of structure theorems/conjectures is about the signs of the coefficients of various expansions of Thom polynomials. 

\subsection{Schur type positivities}

The reader probably already noticed Schur-positivity in many examples in earlier sections.

\begin{theorem} \cite{pragacz:positivity}
    The expansion of Thom series of contact singularities in Schur polynomials have non-negative coefficients.
\end{theorem}

The proof is based on a geometric observation. Therefore the algebraic forms of Thom series, such as those in Sections \ref{sec:finite1}--\ref{sec:finite3}, must automatically produce Schur-positive expansions. This fact imposes convex geometric constraints on the ingredients of those formulas, such as the q-polynomial of Section ~\ref{sec:finite1}.

Schur-positivity seems to extend to the following versions:
\begin{itemize}
    \item The K theory Thom polynomials of Section \ref{sec:KTp}. Here the rule is alternating sign of the Grothendieck polynomial expansion---both in the ``economical'' version of the formula and the ``stable'' version; see \cite{RSz}.
    \item All known SSM Thom polynomials (cf. Section \ref{sec:ssmTp}) of contact singularities satisfy an alternating sign rule in their Schur expansions with respect to cohomological degree, see eg. \eqref{eq:SSM_example}, as well as~\cite{TPP}.
    \item Much fewer coefficients of much fewer motivic Segre Thom polynomials are known, but they all seem to satisfy obvious sign rules in their Grothendieck expansions, see the different examples: \eqref{MSTP_example} and \eqref{eq:MSTP_another}. In analogous situations, besides sign rules, many log-concavity properties have been observed too. Hence it is natural to predict that log-concavity will hold in one way or another for Grothendieck expansions of motivic Segre Thom polynomials.
\end{itemize}

\begin{remark}
Why are the Schur expansion and its K-theory counterpart, the Grothendieck expansion, considered the natural expansions for geometrically relevant polynomials? One answer lies in geometry: Schur and Grothendieck polynomials represent the fundamental classes of Schubert varieties in Grassmannians, which are essential components of geometric structures. However, for SSM Thom polynomials, it may be even more natural to expand in the basis of SSM classes of Schubert varieties (cf. \cite{FRcsm, Sutipoj}). These functions are known as s-tilde functions \cite[\S8.1]{FRcsm}. Indeed, the known SSM Thom polynomials exhibit s-tilde positivity, as demonstrated in \cite{TPP}. The analogous property in K-theory involves MS Thom polynomials expanded in the motivic Segre classes of Schubert varieties.
\end{remark}

\subsection{Monomial positivity for Morin singularities}

In algebraic combinatorics certain polynomials satisfy not only Schur positivity but a stronger property: $e$-positivity. For Thom polynomials this property  translates to positive coefficients, when we expand in Chern monomials. For example 
\[
\begin{array}{ll}
\Tp(A_5,\ell=0)
    & =c_1^5+ 10c_1^3c_2+ 25c_1^2c_3+ 10c_1c_2^2+ 38c_1c_4+ 12c_2c_3+ 24c_5\\
    & =s_{11111}+ 14s_{2111}+ 35s_{221}+ 71s_{311}+ 92s_{32}+ 154s_{41}+ 120s_{5}
\end{array}
\]
is not only Schur positive but also monomial positive. In general, only the weaker Schur positivity holds:
\[
\begin{array}{ll}
\Tp(I_{24},\ell=0)
&=
2c_1^2c_2^2+ 3c_2^3-2c_1^3c_3+ 2c_1c_2c_3-3c_3^2-5c_1^2c_4+ 9c_2c_4 \\
&=2s_{2211}+ 5s_{222}+ 12s_{321}+ 4s_{33}+ 16s_{42}+ 6s_{51}+ 6s_{6}.
\end{array}
\]
As observed in \cite{rrtp}, for singularities that are not of type $A_k$ (i.e., not ``Morin singularities"), the monomial expansion necessarily includes both positive and negative coefficients. This leaves the case of Morin singularities. Interestingly, all known Thom polynomials of Morin singularities are monomial positive \cite[Conj. 5.5]{rrtp}. Although this conjecture has remained open for over 20 years, it has established connections to other areas of mathematics, such as the Green-Griffiths-Lang conjecture and hyperbolicity questions \cite{Bgg}.

\section{What else should have been covered in this paper?}
\label{sec:whatelse}

We conclude this informal survey on Thom polynomials with a list of topics that deserve to be mentioned.

\subsection{Singularities of maps with structures} 
We studied singularities of ``general'' maps. In certain applications, however, the natural maps are not generic; they possess additional structures. For example, they may be Lagrange maps, Legendre maps, or have specific symmetries, or special incidences with another object. We direct the reader to \cite{KazaLL,KazaMulti} and references therein.

\subsection{Thom polynomials of right-left singularities}
We discussed right-left singularities in Section \ref{sec:Asing}, but we did not explore the properties of their Thom polynomials or their relevance in geometry. For more information, see \cite{TO_Asing} and the references therein.

\subsection{The case $\ell<0$} Maps from larger to smaller dimensional spaces have positive-dimensional fibers, and consequently, both their local and global singularity theory behave differently. Thom polynomial theory exists in this context as well; see for example \cite{KazaSapporo}. Matszangosz and Feh\'er announced that they found the analogue of Thom series for singularities with $\ell<0$.

\subsection{Avoiding ideals} \cite{FR04, terpai}
The Thom polynomial is the fundamental class of a singularity locus, and hence it is naturally supported on the singularity locus. There are, however, other polynomials in the same variables that are also supported on the same locus. These polynomials form an ideal known as the ``avoiding ideal'' of the singularity. If \emph{any} element of the avoiding ideal of $\eta$ is non-zero for a map $f: M \to N$, then the map $f$ must contain $\eta$ points (or worse)—meaning that $\eta$ cannot be avoided. An avoiding ideal shares a property with a homogeneous principal ideal: its lowest degree part is one-dimensional, spanned by the Thom polynomial. Determining all the ideal generators of avoiding ideals is an intriguing question with promising applications.

\subsection{Hierarchy of singularities measured by Thom polynomials} \cite{rrtp,zsolt}
Understanding concrete Thom polynomials has implications for classical local singularity theory, too. In particular, Thom polynomials contribute to understanding the hierarchy of singularities. By knowing the Thom polynomial of $\eta$ and the symmetries of $\zeta$, one can define and compute an ``incidence class" $I(\eta, \zeta)$. If this class is non-zero, it indicates that $\eta$ singularities must exist in any neighborhood of $\zeta$, signifying their adjacency in the hierarchy of singularities.

\subsection{Multi-singularity formulas}
The theory of Thom polynomials can be extended to multi-singularities. This extension is arguably the most powerful interface of Thom polynomials with enumerative geometry. A multi-singularity is a multiset of singularities $\eta = \{\eta_1, \eta_2, \ldots, \eta_r\}$ with a distinguished element $\eta_1$. The $\eta$ points of a map $f: M^m \to N^n$ are the points $x_1$ in $M$ such that $f(x_1)$ has exactly $r$ preimages $x_1, \ldots, x_r$, and the singularity of $f$ at $x_i$ is $\eta_i$. In multi-singularity theory, one seeks universal formulas for multi-singularity loci and their $f$-images. The simplest case is $\eta = \{A_0, A_0\}$, where we aim to find the cohomology class of the double point locus of a map $f: M \to N$. The double point formula
\[
[\{A_0, A_0\}(f)] = f^*([f(M)]) - c_{\ell}(f)
\]
is just the tip of the iceberg. Recent developments in multi-singularity formulas reveal connections to interpolation \cite{KazaMulti,RM4tuple}, tautological integrals on Hilbert schemes \cite{BSzmult}, and algebraic cobordism \cite{TOmult}.

\bibliographystyle{alpha}
\bibliography{sample}

\begin{thebibliography}{AGLV98}

\bibitem[AF23]{AF}
D.~Anderson and W.~Fulton.
\newblock {\em Equivariant cohomology in algebraic geometry}.
\newblock CUP, 2023.

\bibitem[AGLV98]{AVGL}
V.~I. Arnold, V.~V. Goryunov, O.~V. Lyashko, and V.~A. Vasil'ev.
\newblock {\em Singularity theory I}.
\newblock Springer, 1993, 1998.

\bibitem[Alu13]{aluffi}
P.~Aluffi.
\newblock Euler characteristics of general linear sections and polynomial {C}hern classes.
\newblock {\em Rend. Circ. Mat. Palermo (special issue)}, pages 3--26, 2013.

\bibitem[B{\'e}r]{GBnr}
G.~B{\'e}rczi.
\newblock Non-reductive geometric invariant theory and {T}hom polynomials.
\newblock arXiv:\-2012.\-06425.

\bibitem[B{\'e}r19]{Bgg}
G.~B{\'e}rczi.
\newblock Thom polynomials and the {G}reen–{G}riffiths–{L}ang conjecture for hypersurfaces with polynomial degree.
\newblock {\em International Mathematics Research Notices}, 2019(22):7037--7092, 2019.

\bibitem[BFR03]{a3}
G.~B\'erczi, L.~Feh\'er, and R.~Rim\'anyi.
\newblock Expressions for resultants coming from the global theory of singularities.
\newblock In L.~McEwan, J.-P. Brasselet, C.~Melles, G.~Kennedy, and K.~Lautier, editors, {\em Topics in Algebraic and Noncommutative Geometry}, number 324 in Contemporary Mathematics. AMS, 2003.

\bibitem[BH61]{borel-haefliger}
A.~Borel and A.~Haefliger.
\newblock La classe d'homologie fondamentale d'un espace analytique.
\newblock {\em Bull. Soc. Math. France}, 89:461--513, 1961.

\bibitem[BS12]{BSz}
G.~B\'erczi and A.~{Sz}enes.
\newblock Thom polynomials of {M}orin singularities.
\newblock {\em Annals of Math.}, 175:567--–629, 2012.

\bibitem[BS21]{BSzmult}
G.~B{\'e}rczi and A.~Szenes.
\newblock Multiple-point residue formulas for holomorphic maps.
\newblock arXiv:\-2112\-.15502, 2021.

\bibitem[BSY10]{BSY}
J.-P. Brasselet, J.~Sch{\"u}rmann, and S.~Yokura.
\newblock Hirzebruch classes and motivic {C}hern classes for singular spaces.
\newblock {\em J. Topol. Anal.}, 2(1):1--55, 2010.

\bibitem[CST]{CsSzT}
A.~Cs\'epai, A.~Sz{\H u}cs, and T.~Terpai.
\newblock On coincidences of singular loci of types {$\Sigma^{1_r}$} and {$\Sigma^r$}.
\newblock Preprint, arXiv:2403.00332.

\bibitem[Dam72]{damon}
J.~Damon.
\newblock Thom polynomials for contact class singularities.
\newblock Ph.D. Thesis, Harvard, 1972.

\bibitem[dPW95]{dPW}
A.~du~Plessis and C.~T.~C. Wall.
\newblock {\em The geometry of topological stability}.
\newblock Oxford Univ. Press, 1995.

\bibitem[EG98]{EG98}
D.~Edidin and W.~Graham.
\newblock Equivariant intersection theory.
\newblock {\em Inv. Math.}, 131(3):595--634, 1998.

\bibitem[Feh21]{LFmc}
L.~M. Feh{\'e}r.
\newblock Motivic {C}hern classes of cones.
\newblock In Javier Fern{\'a}ndez~de Bobadilla, Tam{\'a}s L{\'a}szl{\'o}, and Andr{\'a}s Stipsicz, editors, {\em Singularities and Their Interaction with Geometry and Low Dimensional Topology}, pages 181--205. Springer, 2021.

\bibitem[FM24]{FMrealcsm}
L.~M. Feh\'er and A.~Matszangosz.
\newblock Obstructions for {M}orin and fold maps: {S}tiefel-{W}hitney classes and {E}uler characteristics of singularity loci.
\newblock In preparation, 2024.

\bibitem[FP09]{zsolt}
L.~M. Fehér and Zs. Patakfalvi.
\newblock The incidence class and the hierarchy of orbits.
\newblock {\em Open Mathematics}, 7(3):429--441, 2009.

\bibitem[FR02]{FRint}
L.~Feh{\'e}r and R.~Rim{\'a}nyi.
\newblock Thom polynomials with integer coefficients.
\newblock {\em Illinois J. Math}, 46(4):1145--1158, 2002.

\bibitem[FR04]{FR04}
L.~M. Feh{\'e}r and R.~Rim{\'a}nyi.
\newblock Calculation of {T}hom polynomials and other cohomological obstructions for group actions.
\newblock In {\em Real and Complex Singularities (Sao Carlos, 2002)}, Contemp. Math., pages 69--93. AMS, 2004.

\bibitem[FR07]{dstab}
L.~M. Feh{\'e}r and R.~Rim{\'a}nyi.
\newblock On the structure of {T}hom polynomials of singularities.
\newblock {\em Bull. London Math. Soc.}, 39:541--549, 2007.

\bibitem[FR12]{FRannals}
L.~M. Feh{\'e}r and R.~Rim{\'a}nyi.
\newblock Thom series of contact singularities.
\newblock {\em Annals of Math.}, 176(3):1381--1426, 2012.

\bibitem[FR18]{FRcsm}
L.~M. Feh\'er and R.~Rim\'anyi.
\newblock Chern-{S}chwartz-{M}ac{P}herson classes of degeneracy loci.
\newblock {\em Geometry and Topology}, 22:3575--3622, 2018.

\bibitem[FR24]{LFRRmsThom}
L.~M. Feh\'er and R.~Rim{\'a}nyi.
\newblock Motivic {S}egre {T}hom polynomials.
\newblock In preparation, 2024.

\bibitem[FRW21]{FRWmc}
L.~M. Feh\'er, R.~Rim\'anyi, and A.~Weber.
\newblock Motivic {C}hern classes and {K}-theoretic stable envelopes.
\newblock {\em Proc. London Math. Soc.}, 122:153--189, 2021.

\bibitem[Kaz]{kaza:noas}
M.~\'E. Kazarian.
\newblock Non-associative {H}ilbert scheme and {T}hom polynomials.
\newblock unpublished.

\bibitem[Kaz97]{Kaz97}
M.~E. Kazarian.
\newblock Characteristic classes of singularity theory.
\newblock In {\em The Arnold-Gelfand mathematical seminars}, pages 325--340. Birkhauser, 1997.

\bibitem[Kaz03a]{KazaLL}
M.~Kazarian.
\newblock Thom polynomials for {L}agrange, {L}egendre, and critical point function singularities.
\newblock {\em Proc. London Math. Soc.}, 86(3):707--734, 2003.

\bibitem[Kaz03b]{KazaMulti}
M.~\'E. Kazarian.
\newblock Multisingularities, cobordisms, and enumerative geometry. ({R}ussian).
\newblock {\em Uspekhi Mat. Nauk}, 4(352):29--88, 2003.
\newblock translation in Russian Math. Surveys 58 (2003), no. 4, 665--724.

\bibitem[Kaz06]{KazaSapporo}
M.~E. Kazarian.
\newblock Thom polynomials.
\newblock In {\em Singularity Theory and Its Applications, (Sapporo, 2003)}, number~43 in Adv. Studies in Pure Math., pages 85--136. Math. Soc. Japan, 2006.

\bibitem[KRW20]{KRWell}
S.~Kumar, R.~Rim\'anyi, and A.~Weber.
\newblock Elliptic classes of {S}chubert varieties.
\newblock {\em Math. Annalen}, 378:703--728, 2020.

\bibitem[MO19]{MO}
D.~Maulik and A.~Okounkov.
\newblock {\em Quantum Groups and Quantum Cohomology}, volume 408 of {\em Ast\'erisque}.
\newblock SMF, 2019.

\bibitem[MR10]{RM4tuple}
R.~Marangell and R.~Rim\'anyi.
\newblock The general quadruple point formula.
\newblock {\em American J. of Math.}, 132(4):867--896, 2010.

\bibitem[MS04]{MS}
E.~Miller and B.~Sturmfels.
\newblock {\em Combinatorial Commutative Algebra}, volume 227 of {\em GTM}.
\newblock Springer, 2004.

\bibitem[NO]{Nekarda}
S.~Nekarda and T.~Ohmoto.
\newblock Computing higher {T}hom polynomials for multi-singularity classes.
\newblock In preparation. 2024.

\bibitem[Ohm16]{toru:SMTP}
T.~Ohmoto.
\newblock Singularities of maps and characteristic classes.
\newblock In {\em School on Real and Complex Singularities in S{\~a}o Carlos, 2012}, number~68 in Adv. Studies in Pure Math., pages 191--265, 2016.

\bibitem[Ohm24a]{OhmotoSurvey}
T.~Ohmoto.
\newblock Thom polynomials for singularities of maps.
\newblock In preparation, 2024.

\bibitem[Ohm24b]{TOmult}
T.~Ohmoto.
\newblock Universal polynomials for multi-singularity loci maps.
\newblock In preparation, 2024.

\bibitem[Por71]{porteous}
I.~Porteous.
\newblock Simple singularities of maps.
\newblock In {\em Proc. Liverpool Singularities I}, number 192 in LNM, pages 286--307. Springer, 1971.

\bibitem[PR22]{Sutipoj}
S.~Promtapan and R.~Rim\'anyi.
\newblock Characteristic classes of symmetric and skew-symmetric degeneracy loci.
\newblock In P.~Aluffi, D.~Anderson, M.~Hering, M.~Mustata, and S.~Payne, editors, {\em Facets of Algebraic Geometry. A Collection in Honor of {W}illiam {F}ulton’s 80th Birthday}, number 472 in LMS Lecture Note Series, pages 254--283. CUP, 2022.

\bibitem[PW07]{pragacz:positivity}
P.~Pragacz and A.~Weber.
\newblock Positivity of {S}chur function expansions of {T}hom polynomials.
\newblock {\em Fundamenta Mathematicae}, 195:85--95, 2007.

\bibitem[Rim00]{rrA4}
R.~Rim\'anyi.
\newblock Computation of the {T}hom polynomial of {$\Sigma^{1111}$} via symmetries of singularities.
\newblock In {\em Real and Complex Singularities}, number 412 in RNM, pages 110--118. Chapman and Hall, 2000.

\bibitem[Rim01]{rrtp}
R.~Rim\'anyi.
\newblock Thom polynomials, symmetries and incidences of singularities.
\newblock {\em Inv. Math.}, 143:499--521, 2001.

\bibitem[Ron71]{ronga}
F.~Ronga.
\newblock Le calcul de la classe de cohomologie enti\`ere duale \`a $\overline{\Sigma}^k$.
\newblock In {\em Proc. of Liverpool Sin\-gu\-la\-rities---Symposium, I (1969/70)}, number 192 in LNM, pages 313--315. Springer, 1971.

\bibitem[RS23]{RSz}
R.~Rim{\'a}nyi and A.~Szenes.
\newblock Residues, {G}rothendieck polynomials, and {K}-theoretic {T}hom polynomials.
\newblock {\em IMRN}, 2023(23):20039--20075, 2023.

\bibitem[RV18]{RV}
R.~Rim{\'a}nyi and A.~Varchenko.
\newblock Equivariant {C}hern-{S}chwartz-{M}ac{P}herson classes in partial flag varieties: interpolation and formulae.
\newblock In {\em Schubert Varieties, Equivariant Cohomology and Characteristic Classes}, IMPANGA2015, pages 225--235. EMS, 2018.

\bibitem[RW20]{RWell}
R.~Rim\'anyi and A.~Weber.
\newblock Elliptic classes of {S}chubert varieties via {B}ott-{S}amelson resolution.
\newblock {\em J.~of Topology}, 13, 2020.

\bibitem[SO18]{TO_Asing}
T.~Sasajima and T.~Ohmoto.
\newblock Thom polynomials in {$\mathcal A$}-classification {I}: counting singular projections of a surface.
\newblock In {\em Schubert Varieties, Equivariant Cohomology and Characteristic Classes}, IMPANGA2015, pages 237--259. EMS, 2018.

\bibitem[Ter09]{terpai}
T.~Terpai.
\newblock Calculation of the avoiding ideal for {$\Sigma^{1,1}$}.
\newblock {\em Banach Center Publications}, 85(1):307--313, 2009.

\bibitem[Tho56]{thomoriginal}
R.~Thom.
\newblock Les singularit{\'e}s des applications diff{\'e}rentiables.
\newblock {\em Ann. Inst. Fourier}, pages 43--87, 1955-56.

\bibitem[TPP]{TPP}
Thom {P}olynomial {P}ortal.
\newblock An online registry of known Thom polynomials, https://\-tpp.\-web.\-unc.edu.

\end{thebibliography}

\end{document}